\documentclass[12pt]{article}
\usepackage{amsmath,amssymb,amsfonts,amsthm}
\usepackage{latexsym}
\usepackage[english]{babel}
\usepackage{showidx}
\usepackage{graphicx}
\usepackage{graphics}
\usepackage{babel,amsmath,amssymb,amsbsy,amsfonts,latexsym,amsthm,mathrsfs}

\def\<{\langle}
\def\>{\rangle}

\def\P{{\mathbb P}}
\def\R{{\mathbb R}}
\def\E{{\mathbb E}}
\def\D{{\mathbb D}}

\def\m #1{m_{#1}}

\numberwithin{equation}{section}
\newtheorem{theorem}{Theorem}[section]

\newtheorem{corollary}[theorem]{Corollary}

\newtheorem{lemma}[theorem]{Lemma}

\newtheorem{hypothesis}[theorem]{Hypothesis}

\newtheorem{proposition}[theorem]{Proposition}
\newtheorem{remark}[theorem]{Remark}

\begin{document}

\title{ \textbf{On the distance between\\
probability density functions}}
\author{ \textsc{Vlad Bally}\thanks{%
Universit\'e Paris-Est, LAMA (UMR CNRS, UPEMLV, UPEC), INRIA, F-77454
Marne-la-Vall\'{e}e, France. Email: \texttt{bally@univ-mlv.fr} }\smallskip \\
\textsc{Lucia Caramellino}\thanks{%
Dipartimento di Matematica, Universit\`a di Roma - Tor Vergata, Via della
Ricerca Scientifica 1, I-00133 Roma, Italy. Email: \texttt{%
caramell@mat.uniroma2.it}}\smallskip\\
}
\date{ }
\maketitle

\begin{abstract}
We give estimates of the distance between the densities of the laws of two
functionals $F$ and $G$ on the Wiener space in terms of the
Malliavin-Sobolev norm of $F-G.$ We actually consider a  more general
framework which allows one to treat with similar (Malliavin type)
methods functionals of a Poisson point measure (solutions of jump type
stochastic equations). We use the above estimates in order to obtain a
criterion which ensures that convergence in distribution implies convergence
in total variation distance; in particular, if the functionals at hand are
absolutely continuous, this implies convergence in $L^{1}$ of the densities.
\end{abstract}

\parindent 0pt



{\textbf{Keywords}}: integration by parts formulas, Riesz transform,
Malliavin calculus, weak convergence, total variation distance.

\medskip

{\textbf{2000 MSC}}: 60H07, 60H30.

\section{Introduction}

In this paper we give estimates of the distance between the densities of the
laws of two functionals $F$ and $G$ on the Wiener space in terms of the
Malliavin-Sobolev norm of $F-G.$ Actually, we consider a slightly more general
framework defined in \cite{bib:[BCl]} or \cite{bib:[BCl2]} which allows one
to treat with similar methods functionals of a Poisson point measure
(solutions of jump type stochastic equations). Such estimates may be used in
order to study the behavior of a diffusion process in short time as it is
done in \cite{bib:[B.C]}. But here we focus on a different application: we
use the above estimates in order to obtain a criterion which guarantees that
convergence in distribution implies convergence in total variation distance;
in particular, if the functionals at hand are absolutely continuous, this
implies convergence in $L^{1}$ of the densities. Moreover, by using some
more general distances, we obtain the convergence of the derivatives of the
density functions as well. The main estimates are given in Theorem \ref%
{MainTh1} in the general framework and in Theorem \ref{MainTh2} in the case
of the Wiener space. The convergence result is given in Theorem \ref{General}
and, for the Wiener space, in Theorem \ref{GeneralW}.

The reader interested in the Wiener space case may go directly to Section %
\ref{wiener}. For functionals on the Wiener space we get one more result
which is in between the Bouleau-Hirsch absolute continuity criterion and the
classical criterion of Malliavin for existence and regularity of the density
of the law of a $d$ dimensional functional $F$: we prove that if $F\in
\D^{2,p}$ with $p>d$ and $\P(\det\sigma _{F}>0)>0$ ($\sigma _{F}$ denoting the Malliavin
covariance matrix of $F$) then, conditionally to $\{\sigma _{F}>0\}$ the law
of $F$ is absolutely continuous and the density is lower semi-continuous. This
regularity property implies that the law of $F$ is locally lower bounded by
the Lebesgue measure and this property turns out to be interesting - see the
joint paper \cite{bib:CLT}.

In the last years number of results concerning the weak convergence of
functionals on the Wiener space using Malliavin calculus and Stein's method
have been obtained by Nourdin, Peccati, Nualart and Poly, see \cite{[NP]},
\cite{[NPy]} and \cite{[NNPy]}. In particular in \cite{[NP]} and \cite%
{[NNPy]} the authors consider functionals living in a finite (and fixed)
direct sum of chaoses and prove that, under a very weak non degeneracy
condition, the convergence in distribution of a sequence of such functionals
implies the convergence in total variation. Our initial motivation was to
obtain similar results for general functionals: we consider a sequence of $d$
dimensional functionals $F_{n},n\in {\mathbb{N}},$ which is bounded in ${%
\mathbb{D}}^{3,p}$ for every $p\geq 1.$ Under a very weak non degeneracy
condition (see (\ref{New4})) we prove that the convergence in distribution
of such a sequence implies the convergence in the total variation distance.
Moreover we prove that if a sequence $F_{n},n\in {\mathbb{N}},$ is bounded
in every ${\mathbb{D}}^{3,p}$, $p\geq 1$, $\lim_{n}F_{n}=F$ in $L^{2}$ and $%
\det \sigma _{F}>0$ a.s., then $\lim_{n}F_{n}=F$ in total variation.
Recently, Malicet and Poly \cite{[MPy]} have proved an alternative version
of this result: if $\lim_{n}F_{n}=F$ in ${\mathbb{D}}^{1,2}$ and $\det
\sigma _{F}>0$ a.s. then the convergence takes place in the total variation
distance.

The paper is organized as follows. In Section \ref{th-abstract}, following
\cite{bib:[BCl]}, we introduce an abstract framework which permits to obtain
integration by parts formulas. In Section \ref{th-distance} we give the main
estimate (the distance between two density functions) in this framework and
in Section \ref{th-convergence} we obtain the convergence results. In
Section \ref{wiener} we come back to the Wiener space framework, so here the
objects and the notations are the standard ones from Malliavin calculus (we
refer to Nualart \cite{bib:[N]} for the general theory). Section \ref%
{th-proof} is devoted to the proof of the main estimate, that is of Theorem %
\ref{MainTh1}. Finally, in Section \ref{jumps} we illustrate our convergence
criterion with an example of jump type equation coming from \cite{bib:[BCl]}.

\section{Main results}

\subsection{Abstract integration by parts framework}

\label{th-abstract}

In this section we briefly recall the construction of integration by parts
formulas for functionals of a finite dimensional noise which mimic the
infinite dimensional Malliavin calculus as done in \cite{bib:[BCl]} and \cite%
{bib:[BCl2]}. We are going to introduce operators that represent the finite
dimensional variant of the derivative and the divergence operators from the
classical Malliavin calculus - and as an outstanding consequence all the
constants which appear in the estimates do not depend on the dimension of
the noise. So, given some constants $c_{i}\in {\mathbb{R}},i=1,...,m$ we
denote by $\mathcal{C}(c_{1},...,c_{m})$ the family of universal constants
which depend on $c_{i},i=1,...,m$\ only. So $C\in \mathcal{C}%
(c_{1},...,c_{m})$ means that $C$ depends on $c_{i},i=1,...,m$\ but on
nothing else in the statement. This is crucial in the following theorems.

\label{sect-diffJ}

On a probability space $(\Omega ,\mathcal{F},{\mathbb{P}})$ we consider a
random variable $V=(V_{1},...,V_{J})$ which represents the basic noise. Here
$J\in {\mathbb{N}}$ is a deterministic integer. For each $i=1,...,J$ we
consider two constants $-\infty \leq a_{i}<b_{i}\leq \infty $ that are
allowed to reach $\infty $. We denote
\begin{equation}\label{Oi}
O_{i}=\{v=(v_{1},...,v_{J}):a_{i}<v_{i}<b_{i}\},\quad i=1,\ldots,J.
\end{equation}
The basic hypothesis is that the law of $V$ is absolutely continuous with
respect to the Lebesgue measure on ${\mathbb{R}}^{J}$ and the density $p_{J}$
is smooth with respect to $v_{i}$ on the set $O_{i}.$ The natural example
which comes on in the standard Malliavin calculus is the Gaussian law on ${%
\mathbb{R}}^{J}$, in which $a_{i}=-\infty $ and $b_{i}=+\infty .$ But we may
also (as an example) take $V_{i}$ independent random variables of
exponential law and here, $a_{i}=0$ and $b_{i}=\infty .$

In order to obtain integration by parts formulas for functionals of $V$, one
performs classical integration by parts with respect to $p_{J}(v)dv.$ But in
order to nullify the border terms in $a_{i} $ and $b_{i}$, it suffices to
take into account suitable ``weights''
\begin{equation*}
\pi _{i}\,:\,{\mathbb{R}}^{J}\rightarrow [0,1],\quad i=1,...,J.
\end{equation*}
We give the precise statement of the hypothesis. But let us first set up the
notations we are going to use. We set $C^{k }({\mathbb{R}}^{d})$ the space
of the functions which are continuously differentiable up to order $k$ and $%
C^{\infty}({\mathbb{R}}^{d})$ for functions which are infinitely
differentiable. We use the subscripts $p$, resp. $b$, to denote functions
having polynomial growth, resp. bounded, together with their derivatives,
and this gives $C^{k }_p({\mathbb{R}}^{d})$, $C^{\infty }_p({\mathbb{R}}%
^{d}) $, $C^{k }_b({\mathbb{R}}^{d})$ and $C^{\infty }_b({\mathbb{R}}^{d})$.
For $k\in{\mathbb{N}}$ and for a multi index $\alpha =(\alpha
_{1},...,\alpha _{k})\in \{1,...,d\}^{k}$ we denote $\vert \alpha \vert =k$
and $\partial _{\alpha }f=\partial _{x^{\alpha _{1}}}...\partial _{x^{\alpha
_{k}}}f$. The case $k=0$ is allowed and gives $\partial_\alpha f=f$. We also
set ${\mathbb{N}}^*={\mathbb{N}}\setminus\{0\}$.

So, throughout this paper, we assume the following assumption does hold.

\medskip

\textbf{Assumption.} \emph{The law of the vector $V=(V_{1},...,V_{J})$ is
absolutely continuous with respect to the Lebesgue measure on $\mathbb{R}%
^{J} $ and we denote with $p_{J}$ the density; we assume that $p_{J}$ has
polynomial growth. We also assume that}

\begin{itemize}
\item[(H0)] \emph{for all $i\in \{1,\ldots,J\}$, $0\leq \pi _{i}\leq 1$, $\pi _{i}\in C^{\infty}_b$ and there exist $-\infty\leq a_i<b_i\leq +\infty$ such that, with $O_i$ defined in (\ref{Oi}),  $\{\pi _{i}>0\}\subset O_{i}$; }

\item[(H1)] \emph{the set $\{p_{J}>0\}$ is open in ${\mathbb{R}}^{J}$ and on
$\{p_{J}>0\}$ we have $\ln p_{J}\in C^\infty$. }
\end{itemize}

We define now the functional spaces and the differential operators.

\medskip

$\square $ \textbf{Simple functionals.} A random variable $F$ is called a
simple functional if there exists $f\in C_{p}^{\infty }({\mathbb{R}}^{J})$
such that $F=f(V)$. We denote through $\mathcal{S}$ the set of simple
functionals.

\medskip

$\square $ \textbf{Simple processes.} A simple process is a random variable $%
U=(U_1,\ldots,U_J)$ in ${\mathbb{R}}^J$ such that $U_{i}\in \mathcal{S}$ for
each $i\in \{1,\ldots,J\}$. We denote by $\mathcal{P}$ the space of the
simple processes. On $\mathcal{P}$ we define the scalar product
\begin{equation*}
\<\cdot,\cdot\>:\mathcal{P}\times\mathcal{P}\to \mathcal{S},\quad
(U,V)\mapsto \left\langle U,V\right\rangle_J =\sum_{i=1}^{J}U_{i}V_{i}.
\end{equation*}

\medskip

$\square $ \textbf{The derivative operator.} We define $D:\mathcal{S}%
\rightarrow \mathcal{P}$ by
\begin{equation}
DF:=(D_{i}F)_{i=1,...,J}\in \mathcal{P}\quad \mbox{where }D_{i}F:=\pi
_{i}(V)\partial _{i}f(V).  \label{f1}
\end{equation}

$\square $ \textbf{The divergence operator.} Let $U=(U_1,\ldots, U_J)\in
\mathcal{P}$, so that $U_{i}\in \mathcal{S}$ and $U_{i}=u_{i}(V)$, for some $%
u_{i}\in C_{p}^{\infty }({\mathbb{R}}^{J})$, $i=1,\ldots,J$. We define $%
\delta \, :\,\mathcal{P}\rightarrow \mathcal{S}$ by
\begin{equation}
\delta (U)=\sum_{i=1}^{J}\delta _{i}(U),\quad \mbox{with}\ \delta
_{i}(U):=-(\partial _{v_{i}}(\pi _{i}u_{i})+\pi _{i}u_{i}1_{O_{i}}\partial
_{v_{i}}\ln p_{J})(V),\ i=1,\ldots,J.  \label{Skorohod}
\end{equation}

Clearly, both $D$ and $\delta$ depend on $\pi $ so a correct notation should
be $D^{\pi }$ and $\delta^\pi$. Since here the weights $\pi _{i}$ are fixed,
we do not mention them in the notation.

\medskip

$\square $ \textbf{The Malliavin covariance matrix.} For $F\in\mathcal{S}^d$%
, the Malliavin covariance matrix of $F$ is defined by
\begin{equation*}
\sigma_F^{k,k^{\prime }}=\langle DF^{k},DF^{k^{\prime }}\rangle
_{J}=\sum_{j=1}^{J}D_{j}F^{k}D_{j}F^{k^{\prime }},\quad k,k^{\prime
}=1,\ldots,d.
\end{equation*}%
We also denote
\begin{equation*}
\gamma_F(\omega )=\sigma ^{-1}_F(\omega ),\quad \omega \in
\{\det\sigma_F>0\}.
\end{equation*}

\medskip

$\square $ \textbf{The Ornstein Uhlenbeck operator.} We define $L:\mathcal{%
S\rightarrow S}$ by
\begin{equation}
L(F)=\delta (DF).  \label{OU}
\end{equation}

$\square $ \textbf{Higher order derivatives and norms. }Let $\alpha =(\alpha
_{1},\ldots ,\alpha _{k})$ be a multi-index, with $\alpha _{i}\in \{1,\ldots
,J\}$, for $i=1,\ldots ,k$ and $|\alpha |=k$. For $F\in \mathcal{S}$, we
define recursively
\begin{equation}
D_{(\alpha _{1},\ldots ,\alpha _{k})}F=D_{\alpha _{k}}(D_{(\alpha
_{1},\ldots ,\alpha _{k-1})}F)\quad \mbox{and}\quad D^{(k)}F=\left(
D_{(\alpha _{1},\ldots ,\alpha _{k})}F\right) _{\alpha _{i}\in \{1,\ldots
,J\}}.  \label{f2}
\end{equation}%
We set $D^{(0)}F=F$ and we notice that $D^{(1)}F=DF$. Remark that $%
D^{(k)}F\in \mathbb{R}^{J\otimes k}$ and consequently we define the norm of $%
D^{(k)}F$ as
\begin{equation}
|D^{(k)}F|=\Big(\sum_{\alpha _{1},\ldots ,\alpha _{k}=1}^{J}|D_{(\alpha
_{1},\ldots ,\alpha _{k})}F|^{2}\Big)^{1/2}.  \label{f3}
\end{equation}%
Moreover, we introduce the following norms for simple functionals: for $F\in
\mathcal{S}$ we set
\begin{equation}
|F|_{1,l}=\sum_{k=1}^{l}|D^{(k)}F|=\sum_{k=0}^{l}|D^{(k)}F|,\qquad
|F|_{l}=|F|+|F|_{1,l}  \label{norm1}
\end{equation}%
and for $F=(F^{1},\ldots ,F^{d})\in \mathcal{S}^{d}$, $|F|_{1,l}=%
\sum_{r=1}^{d}|F^{r}|_{1,l}$ and $|F|_{l}=\sum_{r=1}^{d}|F^{r}|_{l}.$
Finally, for $U=(U_{1},\ldots ,U_{J})\in \mathcal{P}$, we set $%
D^{(k)}U=(D^{(k)}U_{1},\ldots ,D^{(k)}U_{J})$ and we define the norm of $%
D^{(k)}U$ as
\begin{equation*}
|D^{(k)}U|=\Big(\sum_{i=1}^{J}|D^{(k)}U_{i}|^{2}\Big)^{1/2}.
\end{equation*}%
We allow the case $k=0$, giving $|U|=\left\langle U,U\right\rangle _{J}^{1/2}
$. Similarly to (\ref{norm1}), we set $|U|_{l}=\sum_{k=0}^{l}|D^{(k)}U|.$

\medskip

$\square $ \textbf{Localization functions. }As it will be clear in the
sequel we need to introduce some localization random variables as follows.
Consider a random variable $\mathbf{\Theta }\in \mathcal{S}$ taking values
on $[0,1]$ and set
\begin{equation*}
d{\mathbb{P}}_{\mathbf{\Theta }}=\mathbf{\Theta }d{\mathbb{P}}.
\end{equation*}%
${\mathbb{P}}_{\mathbf{\Theta }}$ is a non negative measure (but generally
not a probability measure) and we set ${\mathbb{E}}_{\mathbf{\Theta }}$ the
expectation (integral) w.r.t. ${\mathbb{P}}_{\mathbf{\Theta }}$. For $F\in
\mathcal{S}$, we define%
\begin{equation*}
\Vert F\Vert _{p,\mathbf{\Theta }}={\mathbb{E}}_{\mathbf{\Theta }}(|F|^{p})^{1/p}
\end{equation*}%
and
\begin{equation}
\Vert F\Vert _{1,l,p,\mathbf{\Theta }}^{p}={\mathbb{E}}_{\mathbf{\Theta }%
}(|F|_{1,l}^{p})\quad \mbox{and}\quad \Vert F\Vert _{l,p,\mathbf{\Theta }%
}^{p}
=\Vert F\Vert _{p,\mathbf{%
\Theta }}^{p}+\Vert F\Vert _{1,l,p,\mathbf{\Theta }}^{p}.  \label{NORM1}
\end{equation}%
that is $\Vert \cdot \Vert _{p,\mathbf{\Theta }}$ and $\Vert \cdot \Vert
_{l,p,\mathbf{\Theta }}$ are the standard $L^{p}$ respectively Sobolev norms
in Malliavin calculus with ${\mathbb{P}}$ replaced by the localized measure $%
{\mathbb{P}}_{\mathbf{\Theta }}$. Notice also that $\Vert F\Vert _{1,l,p,%
\mathbf{\Theta }}^{p}$ does not take into account the $L^{p}$ norm of $F$
itself but only of the derivatives of $F.$ This is the motivation of
considering this norm.

Since $|F|_{0}=\left\vert F\right\vert $, one has $\Vert F\Vert _{0,p,%
\mathbf{\Theta }}=\Vert F\Vert _{p,\mathbf{\Theta }}$. In the case $\mathbf{%
\Theta }=1$ we come back to the standard notation: $\Vert F\Vert _{p}={%
\mathbb{E}}(|F|^{p})$ and
\begin{equation}
\Vert F\Vert _{1,l,p}^{p}={\mathbb{E}}(|F|_{1,l}^{p})\quad \mbox{and}\quad
\Vert F\Vert _{l,p}^{p}={\mathbb{E}}(|F|^p+|F|_{1,l}^{p}).  \label{NORM2}
\end{equation}%
Notice also that since $\mathbf{\Theta }\leq 1$ we have%
\begin{equation}
\Vert F\Vert _{1,l,p,\mathbf{\Theta }}\leq \Vert F\Vert
_{1,l,p}\quad \mbox{and}\quad \Vert F\Vert _{l,p,\mathbf{\Theta }%
}\leq \Vert F\Vert _{l,p}.  \label{Main0}
\end{equation}%
For $p\in {\mathbb{N}}$ we set
\begin{equation}
\mathrm{m}_{q,p}(\mathbf{\Theta }):=1\vee \Vert \ln \mathbf{\Theta }\Vert
_{1,q,p,\mathbf{\Theta }}.  \label{Mall1}
\end{equation}%
Since $\mathbf{\Theta }>0$ almost surely with respect to ${\mathbb{P}}_{%
\mathbf{\Theta }}$ the above quantity makes sense.

We will work with localization random variables of the following specific
form. For $a>0$, set $\psi _{a}, \phi _{a}:{\mathbb{R}}\rightarrow {\mathbb{R}}_{+}$
as follows:
\begin{equation}
\begin{array}{rcl}
\psi _{a}(x) & = & 1_{|x|\leq a}+\exp \Big(1-\frac{a^{2}}{a^{2}-(|x|-a)^{2}}%
\Big)1_{a<|x|<2a}, \\
\phi _{a}(x) & = & 1_{|x|\geq a}+\exp \Big(1-\frac{a^{2}}{(2|x|-a)^{2}}\Big)%
1_{a/2< |x|< a}.%
\end{array}
\label{Mall10}
\end{equation}%
The function $\psi _{a}$ is suited to localize around zero and $\phi _{a}$
is suited to localize far from zero. Then $\psi _{a},\phi _{a}\in
C_{b}^{\infty }({\mathbb{R}})$, $0\leq \psi _{a}\leq 1$, $0\leq \phi
_{a}\leq 1$ and we have the following property: for every $p,k\in {\mathbb{N}%
}$ there exists a universal constant $C_{k,p}$ such that for every $x\in {%
\mathbb{R}}_{+}$%
\begin{equation}
\psi _{a}(x)\left\vert (\ln \psi _{a})^{(k)}(x)\right\vert ^{p}\leq \frac{%
C_{k,p}}{a^{pk}}\quad \mbox{and}\quad \phi _{a}(x)\left\vert (\ln \phi
_{a})^{(k)}(x)\right\vert ^{p}\leq \frac{C_{k,p}}{a^{pk}}.  \label{loc9}
\end{equation}%
We consider now $\Theta _{i}\in \mathcal{S}$ and $a_{i}>0,i=1,...,l+l^{%
\prime }$ and define
\begin{equation}
\mathbf{\Theta }=\prod_{i=1}^{l}\psi _{a_{i}}(\Theta _{i})\times
\prod_{i=l+1}^{l+l^{\prime }}\phi _{a_{i}}(\Theta _{i}).  \label{Mall10'}
\end{equation}%
As an easy consequence of (\ref{loc9}) we obtain
\begin{equation}
\mathrm{m}_{q,p}(\mathbf{\Theta })\leq 1\vee \Vert \ln \mathbf{\Theta }\Vert
_{1,q,p,\mathbf{\Theta }}\leq 1\vee C_{p,q}\,\sum_{i=1}^{l+l^{\prime }}\frac{%
1}{a_{i}^{q}}\Vert \Theta _{i}\Vert _{1,q,p,\mathbf{\Theta }}.
\label{Mall11}
\end{equation}%
In particular, if $\Vert \Theta _{i}\Vert _{1,q,p}<\infty
,i=1,...,l+l^{\prime }$ then
\begin{equation}
\mathrm{m}_{q,p}(\mathbf{\Theta })\leq 1+C_{p,q}\,\sum_{i=1}^{l+l^{\prime }}%
\frac{1}{a_{i}^{q}}\Vert \Theta _{i}\Vert _{1,q,p}<\infty .  \label{Mall11'}
\end{equation}

Moreover, given some $q\in {\mathbb{N}},p\geq 1$ we denote%
\begin{equation}
\mathrm{U}_{q,p,\mathbf{\Theta }}(F):=\max \{1,{\mathbb{E}}_{\mathbf{\Theta }%
}((\det \sigma _{F})^{-p})(\left\Vert F\right\Vert _{1,q+2,p,\mathbf{\Theta }%
}+\left\Vert LF\right\Vert _{q,p,\mathbf{\Theta }})\}.  \label{Main1}
\end{equation}%
In the case $\mathbf{\Theta }=1$ we have $\mathrm{m}_{q,p}(\mathbf{\Theta }%
)=1$ and
\begin{equation}
\mathrm{U}_{q,p}(F):=\max \{1,{\mathbb{E}}((\det \sigma
_{F})^{-p})(\left\Vert F\right\Vert _{1,q+2,p}+\left\Vert LF\right\Vert
_{q,p})\}.  \label{Main4}
\end{equation}%
Notice that $\mathrm{U}_{q,p,\mathbf{\Theta }}(F)$ and $\mathrm{U}_{q,p}(F)$
do not involve the $L^{p}$ norm of $F$ but only of its derivatives and of $LF.$

We are now able to state the main result in our paper.

\begin{theorem}
\label{MainTh1} Let $q\in {\mathbb{N}}_{\ast }$. We consider the
localization random variable $\mathbf{\Theta }$ defined in (\ref{Mall10'})
and we assume that for every $p\in {\mathbb{N}}$ one has $\Vert \Theta
_{i}\Vert _{q+2,p}<\infty ,i=1,...,l+l^{\prime }.$ In particular $\mathrm{m}%
_{q+2,p}(\mathbf{\Theta })<\infty .$ Let $\mathrm{U}_{q,p,\mathbf{\Theta }%
}(F)$ be as in (\ref{Main1}).

\medskip

\textbf{A}. Let $F\in \mathcal{S}^{d}$ be such that $\mathrm{U}_{q,p,\mathbf{%
\Theta }}(F)<\infty $ for every $p\in {\mathbb{N}}.$ Then under ${\mathbb{P}}%
_{\mathbf{\Theta }}$ the law of $F$ is absolutely continuous with respect to
the Lebesgue measure. We denote by $p_{F,\mathbf{\Theta }}$ its density and
we have $p_{F,\mathbf{\Theta }}\in C^{q-1}({\mathbb{R}}^{d}).$ Moreover
there exist $C,a,b,p\in \mathcal{C}(q,d)$ such that for every $y\in {%
\mathbb{R}}^{d}$ and every multi index $\alpha =(\alpha _{1},...,\alpha
_{k})\in \{1,...,d\}^{k},k\in \{0,...,q\}$ one has%
\begin{equation}
\left\vert \partial _{\alpha }p_{F,\mathbf{\Theta }}(y)\right\vert \leq C%
\mathrm{U}_{q,p,\mathbf{\Theta }}^{a}(F)\times \mathrm{m}_{q+2,p}^{a}(%
\mathbf{\Theta })\times \big({\mathbb{P}}_\mathbf{\Theta }(\left\vert F-y\right\vert < 2 )%
\big)^{b}.  \label{Main2}
\end{equation}%
\textbf{B}. Let $F,G\in \mathcal{S}^{d}$ be such that $\mathrm{U}_{q+1,p,%
\mathbf{\Theta }}(F),\mathrm{U}_{q+1,p,\mathbf{\Theta }}(G)<\infty $ for
every $p\in {\mathbb{N}}$ and let $p_{F,\mathbf{\Theta }}$ and $p_{G,\mathbf{%
\Theta }}$ be the densities of the laws of $F$ respectively of $G$ under ${%
\mathbb{P}}_{\mathbf{\Theta }}.$ There exist $C,a,b,p\in \mathcal{C}(q,d)$
such that for every $y\in {\mathbb{R}}^{d}$ and every multi index $\alpha
=(\alpha _{1},...,\alpha _{k})\in \{1,...,d\}^{k},$ $0\leq k\leq q$ one has
\begin{equation}  \label{Main3}
\begin{array}{rcl}
\left\vert \partial _{\alpha }p_{F,\mathbf{\Theta }}(y)-\partial _{\alpha
}p_{G,\mathbf{\Theta }}(y)\right\vert & \leq & C\mathrm{U}_{q+1,p,\mathbf{%
\Theta }}^{a}(F)\times \mathrm{U}_{q+1,p,\mathbf{\Theta }}^{a}(G)\times
\mathrm{m}_{q+2,p}^{a}(\mathbf{\Theta })\times \smallskip \\
&  & \times \big({\mathbb{P}}_\mathbf{\Theta }(\left\vert F-y\right\vert < 2\vert)+{\mathbb{P}}_\mathbf{\Theta }(\left\vert G-y\right\vert < 2 )\big)^{b}\times\smallskip \\
&  & \times \left( \Vert F-G\Vert _{q+2,p,\mathbf{\Theta }}+\Vert LF-LG\Vert
_{q,p,\mathbf{\Theta }}\right) .%
\end{array}%
\end{equation}
\end{theorem}

\begin{remark}
The above result can be written in the case $\mathbf{\Theta }=1$. Here, $\mathrm{m}%
_{q+2,p}(\mathbf{\Theta })=1$ and the quantities $\Vert F-G\Vert _{q+2,p,%
\mathbf{\Theta }}$ and $\Vert LF-LG\Vert _{q,p,\mathbf{\Theta }}$ are
replaced by $\Vert F-G\Vert _{q+2,p}$ and $\Vert LF-LG\Vert _{q,p}$
respectively.
\end{remark}

\begin{remark}
Since $\P_\mathbf{\Theta }(A)\leq \P(A)$, in (\ref{Main2}) and (\ref{Main3}) one can replace $\P_\mathbf{\Theta}$ with $\P$.
\end{remark}

\begin{remark}
\label{use-main} Estimates (\ref{Main2}) and (\ref{Main3}) may be rewritten
in terms of the queues of the law of $F$ and $G$ by noticing that if $|y|>4$
then $\{|F-y|< 2\}\subset \{|F|>|y|/2\}$ and $\{|G-y|< 2\}\subset
\{|G|>|y|/2\}$. But we can do something else. In fact, by using the Markov
inequality, for every $\ell\geq 1$ and for $|y|>4$ we get ${\mathbb{P}}_{\mathbf{\Theta}}
(|F-y|<2)\leq {\mathbb{P}}_{\mathbf{\Theta}}(|F|>|y|/2)\leq C\,{\mathbb{E}}_{\mathbf{\Theta}}(|F|^\ell)/(1+|y|)^%
\ell,$ $C$ denoting a universal constant. And by taking into account also
the case $|y|\leq 4$, for a suitable $C$ we have
\begin{equation*}
{\mathbb{P}}_{\mathbf{\Theta}}(|F-y|<2)\leq C\,\frac{1\vee {\mathbb{E}}_{\mathbf{\Theta}}(|F|^\ell)}{(1+|y|)^\ell%
}
\end{equation*}
and similarly for $G$. Then, the second factors in formulas (\ref{Main2})
and (\ref{Main3}) may be written in terms of the above inequality as
follows: for every $\ell\geq 1$ and $y\in{\mathbb{R}}^d$,
\begin{equation}  \label{Main2bis}
\left\vert \partial _{\alpha }p_{F,\mathbf{\Theta }}(y)\right\vert \leq C%
\mathrm{U}_{q,p,\mathbf{\Theta }}^{a}(F)\times \mathrm{m}_{q+2,p}^{a}(%
\mathbf{\Theta })\times \displaystyle\frac{(1+\|F\|_{\ell, {\mathbf{\Theta}}}^\ell)^b}{%
(1+|y|)^{\ell b}}
\end{equation}%
and
\begin{equation}  \label{Main3bis}
\begin{array}{rcl}
\left\vert \partial _{\alpha }p_{F,\mathbf{\Theta }}(y)-\partial _{\alpha
}p_{G,\mathbf{\Theta }}(y)\right\vert & \leq & C\mathrm{U}_{q+1,p,\mathbf{%
\Theta }}^{a}(F)\times \mathrm{U}_{q+1,p,\mathbf{\Theta }}^{a}(G)\times
\mathrm{m}_{q+2,p}^{a}(\mathbf{\Theta })\times \smallskip \\
&  & \times \displaystyle\frac{(1+\|F\|_{\ell, {\mathbf{\Theta}}}^\ell+\|G\|_{\ell, {\mathbf{\Theta}}}^\ell)^b}{%
(1+|y|)^{\ell b}}\times \smallskip \\
&  & \times \left( \Vert F-G\Vert _{q+2,p,\mathbf{\Theta }}+\Vert LF-LG\Vert
_{q,p,\mathbf{\Theta }}\right) .%
\end{array}%
\end{equation}
\end{remark}

The proof of Theorem \ref{MainTh1} is the main effort in our paper and it is
postponed for Section \ref{th-proof} (see Proposition \ref{1} \textbf{C.}
and Theorem \ref{Diference}).

As a consequence of Theorem \ref{MainTh1} we obtain the following
regularization result. Let $\gamma _{\delta }$ be the density of the centred
normal law of covariance $\delta \times I$ on ${\mathbb{R}}^{d}.$ Here $%
\delta >0$ and $I$ is the identity matrix.

\begin{lemma}
\label{Regularization}There exist some universal constants $C,p,a\in
\mathcal{C}(d)$ such that for every $\varepsilon >0,\delta >0$ and every $%
F\in \mathcal{S}^{d}$ one has%
\begin{equation}
\left\vert {\mathbb{E}}(f(F))-{\mathbb{E}}(f\ast \gamma _{\delta
}(F))\right\vert \leq C\left\Vert f\right\Vert _{\infty }\Big({\mathbb{P}}%
(\sigma _{F}<\varepsilon )+\frac{\sqrt{\delta }}{\varepsilon ^{p}}%
(1+\left\Vert F\right\Vert _{3,p}+\left\Vert LF\right\Vert _{1,p})^{a}\Big)
\label{Main5}
\end{equation}%
for every bounded and measurable $f:{\mathbb{R}}^{d}\rightarrow {\mathbb{R}}.
$ Moreover, for if $f\in L^{1}(\R^{d})$%
\begin{equation}
\left\vert {\mathbb{E}}(f(F))-{\mathbb{E}}(f\ast \gamma _{\delta
}(F))\right\vert \leq C(\left\Vert f\right\Vert _{\infty }+\left\Vert
f\right\Vert _{1})\Big({\mathbb{P}}(\sigma _{F}<\varepsilon )+\frac{\sqrt{%
\delta }}{\varepsilon ^{p}}(1+\left\Vert F\right\Vert _{1,3,p}+\left\Vert
LF\right\Vert _{1,p})^{a}\Big)  \label{Main5'}
\end{equation}
\end{lemma}

Notice that in the r.h.s. of (\ref{Main5'}) $\left\Vert F\right\Vert _{3,p}$
is replaced by $\left\Vert F\right\Vert _{1,3,p}$ so $\left\Vert
F\right\Vert _{p}$ is not involved. The price to be payed is that we have to
replace $\left\Vert f\right\Vert _{\infty }$
with  $\left\Vert f\right\Vert _{\infty }+\left\Vert f\right\Vert _{1}$.

\bigskip

\textbf{Proof}. Along this proof $C$\ denotes a constant in $\mathcal{C}(d)$
which may change from a line to another. We construct the localization
random variable $\mathbf{\Theta }_{\varepsilon }=\phi _{\varepsilon }(\det
\sigma _{F})$ with $\phi _{\varepsilon }$ given in (\ref{Mall10}). By (\ref%
{Mall11}) for every $p\geq 1$%
\begin{equation}
\mathrm{m}_{q,p}(\mathbf{\Theta }_{\varepsilon })\leq \frac{C}{\varepsilon
^{q}}\Vert \det \sigma _{F}\Vert _{q,p,\mathbf{\Theta }_{\varepsilon }}\leq
\,\frac{C}{\varepsilon ^{q}}\Vert F\Vert _{1,q+1,p}^{d}.  \label{Main6}
\end{equation}%
We fix $\delta \in (0,1)$ and we define $F_{\delta }=F+\delta \Delta $ where
$\Delta $ is a standard Gaussian random variable independent of $V.$ We will
use the result in Theorem \ref{MainTh1}, here not with respect to $%
V=(V_{1},...,V_{J})$ but with respect to $(V,\Delta
)=(V_{1},...,V_{J},\Delta ).$ The Malliavin covariance matrix of $F$ with
respect to $(V,\Delta )$ is the same as the one with respect to $V$ (because
$F$ does not depend on $\Delta )$ so on the set $\{\mathbf{\Theta }%
_{\varepsilon }\neq 0\}$ we have $\det \sigma _{F}\geq \varepsilon .$ We
denote by $\sigma _{F_{\delta }}$ the Malliavin covariance matrix of $%
F_{\delta }$ computed with respect to $(V,\Delta ).$ We have $\left\langle
\sigma _{F_{\delta }}\xi ,\xi \right\rangle =\delta \left\vert \xi
\right\vert ^{2}+\left\langle \sigma _{F}\xi ,\xi \right\rangle .$ By Lemma
7-29, pg 92 in \cite{bib:BGJ}, for every symmetric non negative defined
matrix $Q$ one has
\begin{equation*}
\frac{1}{\det Q}\leq C_{1}\int_{{\mathbb{R}}^{d}}\left\vert \xi \right\vert
^{d}e^{-\left\langle Q\xi ,\xi \right\rangle }d\xi \leq C_{2}\frac{1}{\det Q}
\end{equation*}%
where $C_{1}$ and $C_{2}$ are universal constants. Using these two
inequalities we obtain $\det \sigma _{F_{\delta }}\geq \frac{1}{C}\det
\sigma _{F}\geq \frac{1}{C}\varepsilon $ on the set $\mathbf{\Theta }%
_{\varepsilon }>0.$ So for $\varepsilon \in (0,1)$ we have
\begin{equation*}
\Vert (\det \sigma _{F})^{-1}\Vert _{p,\mathbf{\Theta }_{\varepsilon
}}+\Vert (\det \sigma _{F_{\delta }})^{-1}\Vert _{p,\mathbf{\Theta }%
_{\varepsilon }}\leq C\varepsilon ^{-1}.
\end{equation*}%
It is also easy to check that
\begin{equation*}
\left\Vert F_{\delta }\right\Vert _{1,3,p,\mathbf{\Theta }_{\varepsilon
}}+\left\Vert LF_{\delta }\right\Vert _{1,p,\mathbf{\Theta }_{\varepsilon
}}\leq C(1+\left\Vert F\right\Vert _{1,3,p,\mathbf{\Theta }_{\varepsilon
}}+\left\Vert LF\right\Vert _{1,p,\mathbf{\Theta }_{\varepsilon }})
\end{equation*}%
so finally we obtain%
\begin{eqnarray*}
\mathrm{U}_{1,p,,\mathbf{\Theta }_{\varepsilon }}(F)+\mathrm{U}_{1,p,,%
\mathbf{\Theta }_{\varepsilon }}(F_{\delta }) &\leq &C(1+\varepsilon
^{-p}(\left\Vert F\right\Vert _{1,3,p,\mathbf{\Theta }_{\varepsilon
}}+\left\Vert LF\right\Vert _{1,p,\mathbf{\Theta }_{\varepsilon }})) \\
&\leq &C(1+\varepsilon ^{-p}(\left\Vert F\right\Vert _{1,3,p}+\left\Vert
LF\right\Vert _{1,p})).
\end{eqnarray*}

By using (\ref{Main6}), we apply Theorem \ref{MainTh1} and we obtain
\begin{eqnarray}
\left\vert p_{F,\mathbf{\Theta }_{\varepsilon }}(y)-p_{F_{\delta },\mathbf{%
\Theta }_{\varepsilon }}(y)\right\vert  &\leq &C(1+\varepsilon
^{-p}(\left\Vert F\right\Vert _{1,3,p}+\left\Vert LF\right\Vert
_{1,p}))^{a}\times   \label{a} \\
&&\times \left( \Vert F-F_{\delta }\Vert _{2,p}+\Vert LF-LF_{\delta }\Vert
_{0,p}\right) .  \notag
\end{eqnarray}%
The r.h.s. of the above inequality does not depend on $y$, so its
integral over $\R^{d}$ is infinite. In order to obtain a finite
integral we use inequality (\ref{Main3bis}) discussed in Remark \ref%
{use-main} with $\ell $ large enough: we may find $C,p,a,b\in \mathcal{C}(d)$
such that
\begin{eqnarray}
\left\vert p_{F,\mathbf{\Theta }_{\varepsilon }}(y)-p_{F_{\delta },\mathbf{%
\Theta }_{\varepsilon }}(y)\right\vert  &\leq &C(1+\varepsilon
^{-p}(\left\Vert F\right\Vert _{3,p}+\left\Vert LF\right\Vert
_{1,p}))^{a}\times \frac{1}{(1+|y|)^{2d}}\times   \label{a1} \\
&&\times \left( \Vert F-F_{\delta }\Vert _{2,p}+\Vert LF-LF_{\delta }\Vert
_{0,p}\right) .  \notag
\end{eqnarray}%
But now $\left\Vert F\right\Vert _{p}$ comes on and this is why we have to
replace $\left\Vert F\right\Vert _{1,3,p}$ by $\left\Vert F\right\Vert
_{3,p}.$

Moreover one can easily check using directly the definitions that
\begin{equation*}
\Vert F-F_{\delta }\Vert _{2,p}+\Vert LF-LF_{\delta }\Vert _{0,p}\leq
C\delta ^{1/2}.
\end{equation*}%
So finally we obtain%
\begin{equation}
\left\vert p_{F,\mathbf{\Theta }_{\varepsilon }}(y)-p_{F_{\delta },\mathbf{%
\Theta }_{\varepsilon }}(y)\right\vert \leq \frac{C}{(1+\left\vert
y\right\vert )^{2d}\varepsilon ^{p}}(1+\left\Vert F\right\Vert
_{3,p}+\left\Vert LF\right\Vert _{1,p}))^{a}\times \sqrt{\delta }.
\label{Main7}
\end{equation}

We are now ready to start the proof of our Lemma. We take $f\in C({\mathbb{R}%
}^{d})$ with $\left\Vert f\right\Vert _{\infty }<\infty$ and we write%
\begin{align*}
{\mathbb{E}}(f(F))-{\mathbb{E}}(f\ast \gamma _{\delta }(F))
=&{\mathbb{E}}%
(f(F))-{\mathbb{E}}(f(F_{\delta })) \\
=&\big[{\mathbb{E}}(f(F)(1-\mathbf{\Theta }_{\varepsilon }))-{\mathbb{E}}%
(f(F_{\delta })(1-\mathbf{\Theta }_{\varepsilon }))\big]+ \\
&+\big[{\mathbb{E}}(f(F)\mathbf{\Theta }_{\varepsilon })-{\mathbb{E}}%
(f(F_{\delta })\mathbf{\Theta }_{\varepsilon })\big] \\
=:&I(\delta ,\varepsilon )+J(\delta ,\varepsilon ).
\end{align*}%
We have%
\begin{equation*}
\left\vert I(\delta ,\varepsilon )\right\vert \leq 2\left\Vert f\right\Vert _{\infty }{\mathbb{E}}(\left\vert 1-%
\mathbf{\Theta }_{\varepsilon }\right\vert) \leq 2\left\Vert f\right\Vert _{\infty }{\mathbb{P}}(\det \sigma
_{F}<\varepsilon ).
\end{equation*}%
We use (\ref{Main7}) in order to obtain%
\begin{eqnarray*}
\left\vert J(\delta ,\varepsilon )\right\vert  &=&\left\vert {\mathbb{E}}_{%
\mathbf{\Theta }_{\varepsilon }}(f(F))-{\mathbb{E}}_{\mathbf{\Theta }%
_{\varepsilon }}(f(F_{\delta }))\right\vert  \\
&=&\left\vert \int f(x)(p_{F,\mathbf{\Theta }_{\varepsilon
}}(x)-p_{F_{\delta },\mathbf{\Theta }_{\varepsilon }}(x))dx\right\vert \leq
\left\Vert f\right\Vert _{\infty }\int \left\vert p_{F,\mathbf{\Theta }_{\varepsilon }}(x)-p_{F_{\delta },%
\mathbf{\Theta }_{\varepsilon }}(x)\right\vert dx \\
&\leq &\frac{C}{\varepsilon ^{p}}\,\left\Vert f\right\Vert _{\infty }(1+\left\Vert F\right\Vert
_{3,p}+\left\Vert LF\right\Vert _{1,p}))^{a}\times \sqrt{\delta }\times \int
\frac{1}{(1+\left\vert y\right\vert )^{2d}}dy
\end{eqnarray*}%
and (\ref{Main5}) follows. We write now
\begin{eqnarray*}
\left\vert J(\delta ,\varepsilon )\right\vert  &=&\left\vert \int f(x)(p_{F,%
\mathbf{\Theta }_{\varepsilon }}(x)-p_{F_{\delta },\mathbf{\Theta }%
_{\varepsilon }}(x))dx\right\vert  \\
&\leq &\left\Vert p_{F,\mathbf{\Theta }_{\varepsilon }}-p_{F_{\delta },%
\mathbf{\Theta }_{\varepsilon }}\right\Vert _{\infty }\left\Vert
f\right\Vert _{1}.
\end{eqnarray*}%
Using (\ref{a}) we obtain (\ref{Main5'}). $\square $

\subsection{Distances and basic estimate}

\label{th-distance}

In this section we discuss the convergence in the total variation distance
defined by%
\begin{equation*}
d_{TV}(F,G)=\sup \{\left\vert {\mathbb{E}}(f(F))-{\mathbb{E}}%
(f(G))\right\vert :\left\Vert f\right\Vert _{\infty }\leq 1\}
\end{equation*}%
The convergence in this distance is related to the convergence of the
densities of the laws: given a sequence of random variables $F_{n}\sim
p_{n}(x)dx$ and $F\sim p(x)dx$ then $d_{TV}(F_{n},F)\rightarrow 0$ is
equivalent to
\begin{equation*}
\lim_{n}\int \left\vert p(x)-p_{n}(x)\right\vert dx=0.
\end{equation*}%
We also consider the Fortet-Mourier distance defined by%
\begin{equation*}
d_{FM}(F,G)=\sup \{\left\vert {\mathbb{E}}(f(F))-{\mathbb{E}}%
(f(G))\right\vert :\left\Vert f\right\Vert _{\infty }+\left\Vert \nabla
f\right\Vert _{\infty }\leq 1\}
\end{equation*}%
and the Wasserstein distance%
\begin{equation*}
d_{W}(F,G)=\sup \{\left\vert {\mathbb{E}}(f(F))-{\mathbb{E}}%
(f(G))\right\vert :\left\Vert \nabla f\right\Vert _{\infty }\leq 1\}.
\end{equation*}

The convergence in $d_{W}$ is equivalent to the convergence in distribution
plus the convergence of the first order moments. Clearly $d_{FM}(F,G)\leq
d_{W}(F,G)$ so convergence in distribution plus the convergence of the first
order moments implies convergence in $d_{FM}.$ One also has $d_{FM}(F,G)\leq
d_{TV}(F,G).$ The aim of this section is to prove a kind of converse type
inequality.

We will be interested in a larger class of distances that we define now. For
$f\in C_{b}^{m}({\mathbb{R}}^{d})$ we denote
$$\left\Vert f\right\Vert
_{m,\infty }=\left\Vert f\right\Vert _{\infty }+\sum_{1\leq \left\vert
\alpha \right\vert \leq m}\left\Vert \partial _{\alpha }f\right\Vert
_{\infty }.$$
Then we define%
\begin{equation}
d_{m}(F,G)=\sup \{\left\vert {\mathbb{E}}(f(F))-{\mathbb{E}}%
(f(G))\right\vert :\left\Vert f\right\Vert _{m,\infty }\leq 1\}.
\label{Main9}
\end{equation}%
So
\begin{equation*}
d_{FM}=d_{1}\qquad \mbox{and}\qquad d_{TV}=d_{0}.
\end{equation*}
Our basic estimate is the following. For $F\in \mathcal{S}^{d}$ we denote%
\begin{equation}
A_{l}(F):=\left\Vert F\right\Vert _{3,l}+\left\Vert LF\right\Vert _{1,l}
\label{Main8}
\end{equation}

\begin{theorem}
\label{Distance}Let $k\in {\mathbb{N}}.$ There exist universal constants $%
C,l,b\in \mathcal{C}(d,k)$ such that for every $F,G\in \mathcal{S}^{d}$ with
$A_{p}(F),A_{p}(G)<\infty ,\forall p\in {\mathbb{N}},$ and every $%
\varepsilon >0$ one has%
\begin{equation}  \label{New1}
\begin{array}{rcl}
d_{0}(F,G) & \leq & \displaystyle\frac{C}{\varepsilon ^{b}}%
(1+A_{l}(F)+A_{l}(G))^{b}d_{k}^{\frac{1}{k+1}}(F,G)+ \smallskip \\
&  & +C{\mathbb{P}}(\det \sigma _{F} <\varepsilon )+C{\mathbb{P}}(\det
\sigma _{G}<\varepsilon ).%
\end{array}%
\end{equation}
\end{theorem}

\textbf{Proof}. Let $\delta >0$ and let $f\in C({\mathbb{R}}^{d})$ with $%
\left\Vert f\right\Vert _{\infty }\leq 1.\ $Since $\left\Vert f\ast \gamma
_{\delta }\right\Vert _{k,\infty }\leq C\delta ^{-k/2}$ we have
\begin{equation*}
\left\vert {\mathbb{E}}(f\ast \gamma _{\delta }(F))-{\mathbb{E}}(f\ast
\gamma _{\delta }(G))\right\vert \leq C\delta ^{-k/2}d_{k}(F,G).
\end{equation*}%
Then using (\ref{Main5})%
\begin{eqnarray*}
\left\vert {\mathbb{E}}(f(F))-{\mathbb{E}}(f(G))\right\vert &\leq &C\delta
^{-k/2}d_{k}(F,G)+C{\mathbb{P}}(\det \sigma _{F}<\varepsilon )+C{\mathbb{P}}%
(\det \sigma _{G}<\varepsilon )+ \\
&&+\frac{C\delta ^{1/2}}{\varepsilon ^{p}}(1+A_{l}(F)+A_{l}(G))^{a}.
\end{eqnarray*}%
We optimize over $\delta $: we take%
\begin{equation*}
\delta ^{(k+1)/2}=d_{k}(F,G) \Big(\frac{1}{\varepsilon^p} (1+
A_{l}(F)+A_{l}(G))^a\Big)^{-1}.
\end{equation*}%
We insert this in the previous inequality and we obtain (\ref{New1}). $%
\square $

\subsection{Convergence results}

\label{th-convergence}

In the previous sections we considered a functional $F\in \mathcal{S}^{d}$
with $\mathcal{S}$ associated to a certain random variable $%
V=(V_{1},...,V_{J}).$ So $F=f(V).$ But the estimates that we have obtained
are estimates of the law and so it is not necessary that the random
variables at hand are functionals of the same $V.$ We may have $F=f(V)$ and $%
\overline{F}=\overline{f}(\overline{V})$ with $\overline{V}=(\overline{V}%
_{1},...,\overline{V}_{\overline{J}}).$ Having this in mind, for a fixed
random variable $V=(V_{1},...,V_{J})$ we denote by $\mathcal{S}%
(V)=\{F=f(V):f\in C_{b}^{\infty }({\mathbb{R}}^{J})\}$ the space of the
simple functionals associated to $V.$ We denote by $\sigma _{F}(V)$ the
Malliavin covariance matrix and
\begin{equation*}
A_{p}(V,F):=\left\Vert F\right\Vert _{3,p}+\left\Vert LF\right\Vert _{1,p}.
\end{equation*}%
Here the norms $\left\Vert F\right\Vert _{q,l}$ and the operator $LF$ are
defined as in (\ref{norm1}) and (\ref{OU}) with respect to $V$.

In the following we will work with a sequence $(F_{n})_{n\in {\mathbb{N}}}$
of $d$ dimensional functionals $F_{n}=(F_{n,1},...,F_{n,d}).$ For each $n$, $%
F_{n,i}\in \mathcal{S}(V_{(n)}),$ $i=1,...,d$ for some random variables $%
V_{(n)}=(V_{(n),1},...,V_{(n),J_{n}}).$ We will use the following two
assumptions. First, we consider a regularity assumption:%
\begin{equation}
\overline{F}_{p}:=\sup_{n}A_{p}(V_{(n)},F_{n})<\infty ,\qquad \forall p\geq
1.  \label{New3}
\end{equation}%
The second one is a (very weak) non degeneracy hypothesis:%
\begin{equation}
\lim_{\varepsilon \rightarrow 0}\eta (\varepsilon )=0\qquad \mbox{with}%
\qquad \eta (\varepsilon ):=\limsup_n{\mathbb{P}}(\det \sigma
_{F_{n}}(V_{(n)})\leq \varepsilon )  \label{New4}
\end{equation}

One has

\begin{lemma}
Let $\overline{F}_p$ be as in (\ref{New3}). If $\overline{F}_1<\infty$ then (%
\ref{New4}) is equivalent to
\begin{equation}
\lim_{\varepsilon \rightarrow 0}\overline{\eta }(\varepsilon )=0\qquad %
\mbox{with}\qquad \overline{\eta }(\varepsilon ):=\limsup_n{\mathbb{P}}%
(\lambda (F_{n})\leq \varepsilon )  \label{New5}
\end{equation}%
where $\lambda (F_{n})$ is the smaller eigenvalue of $\sigma
_{F_{n}}(V_{(n)}).$
\end{lemma}

\textbf{Proof}. The statement is trivial for $d=1$, so we consider the case $%
d>1$. Since $\det \sigma _{F_{n}}(V_{(n)})\geq \lambda ^{d}(F_{n})$ we have $%
{\mathbb{P}}(\det \sigma _{F_{n}}(V_{(n)})\leq \varepsilon )\leq {\mathbb{P}}%
(\lambda (F_{n})\leq \varepsilon ^{1/d})$ so that $\eta (\varepsilon )\leq
\overline{\eta }(\varepsilon ^{1/d})$. If $\gamma (F_{n})$ is the largest
eigenvalue of $\sigma _{F_{n}}(V_{(n)})$ then $\det \sigma
_{F_{n}}(V_{(n)})\leq \lambda (F_{n})\gamma ^{d-1}(F_{n})$ so that%
\begin{align*}
{\mathbb{P}}(\lambda (F_{n}) \leq \varepsilon ) &\leq {\mathbb{P}}(\det
\sigma _{F_{n}}(V_{(n)})\leq \varepsilon \gamma ^{d-1}(F_{n})) \\
&\leq {\mathbb{P}}(\gamma ^{d-1}(F_{n})\geq \varepsilon ^{-1/2})+{\mathbb{P}}%
(\det \sigma _{F_{n}}(V_{(n)})\leq \varepsilon ^{1/2})
\end{align*}
But $\gamma(F_n)\leq |DF_n|^2$, so
\begin{equation*}
{\mathbb{P}}(\gamma ^{d-1}(F_{n})\geq \varepsilon ^{-1/2}) \leq
\varepsilon^{\frac 1{4(d-1)}}{\mathbb{E}}(|DF_n|)\leq \varepsilon^{\frac
1{4(d-1)}}\,\overline{F}_1.
\end{equation*}
We conclude that $\overline{\eta }(\varepsilon )\leq \varepsilon ^{\frac
1{4(d-1)}} \overline{F}_1+\eta (\varepsilon ^{1/2}).$ $\square $

\begin{theorem}
\label{General} We consider a sequence of functionals $%
F_{n}=(F_{n,1},...,F_{n,d})\in \mathcal{S}^{d}(V_{(n)})$ and we assume that  (\ref{New3}) and (\ref{New4}) hold. Suppose also that $\lim_{n}F_{n}=F$ in
distribution and $\lim_{n}\E(F_{n})=\E(F)$. Then
\begin{equation}
\lim_{n}d_{TV}(F,F_{n})=0.  \label{New6}
\end{equation}%
In particular if the laws of $F$ and $F_{n}$ are absolutely continuous with
density $p_{F}$ and $p_{F_{n}}$ respectively, then
\begin{equation*}
\lim_{n}\int \left\vert p_{F}(x)-p_{F_{n}}(x)\right\vert dx=0.
\end{equation*}
\end{theorem}

\textbf{Proof}. Using (\ref{New1}) with $k=1$ we may find some $C,l,b\in
\mathcal{C}(d,k)$ such that for every $n,m\in {\mathbb{N}}$
\begin{equation*}
d_{0}(F_{n},F_{m})\leq \frac{C}{\varepsilon ^{b}}(1+\overline{F}%
_{l})^{b}\,d_{1}^{1/2}(F_{n},F_{m})+C{\mathbb{P}}(\det \sigma
_{F_{n}}(V_{(n)})<\varepsilon )+C{\mathbb{P}}(\det \sigma
_{F_{m}}(V_{(m)})<\varepsilon ). 
\end{equation*}
Since $\lim_{n}F_{n}=F$ in law then $\limsup _{n,m\rightarrow \infty
}d_{1}(F_{n},F_{m})=0$, so that $\limsup _{n,m\rightarrow \infty
}d_{0}(F_{n},F_{m})\leq C\eta (\varepsilon ).$ This is true for every $%
\varepsilon >0.$ So using (\ref{New4}) we obtain $\limsup _{n,m\rightarrow
\infty }d_{0}(F_{n},F_{m})$ $=0.$ $\square $

\subsection{Functionals on the Wiener space}

\label{wiener}

Let $(\Omega ,{\mathscr F},{\mathbb{P}})$ be a probability space where a
Brownian motion $W=(W^{1},...,W^{N})$ is defined. We briefly recall the main
notations in Malliavin calculus, for which we refer to Nualart \cite{bib:[N]}%
. We denote by ${\mathbb{D}}^{m,p}$ the space of the random variables which
are $m$ times differentiable in Malliavin sense in $L^{p}$ and for a
multi-index $\alpha =(\alpha _{1},\ldots ,\alpha _{k})\in \{1,\ldots
,N\}^{k} $, $k\leq m$, we denote by $D^{\alpha }F$ the Malliavin derivative
of $F$ corresponding to the multi-index $\alpha .$ Moreover we define%
\begin{equation}
|D^{(k)}F|^{2}=\sum_{|\alpha |=k}\int_{[0,1)^{k}}\left\vert D_{s_{1},\ldots
,s_{k}}^{\alpha }F\right\vert ^{2}\,ds_{1},\ldots ds_{k}\quad \mbox{and}%
\quad \left\vert F\right\vert _{m}^{2}=\left\vert F\right\vert^2
+\sum_{k=1}^{m}|D^{(k)}F|^{2}.  \label{W4}
\end{equation}%
So, ${\mathbb{D}}^{m,p}$ is the closure of the space of the simple
functionals with respect to the Malliavin Sobolev norm
\begin{equation}
\left\Vert F\right\Vert _{m,p}^{p}={\mathbb{E}}\big(|F|_{m}^{p}\big)
\label{W3}
\end{equation}%
We set ${\mathbb{D}}^{m,\infty }=\cap _{p\geq 1}{\mathbb{D}}^{m,p}$ and ${%
\mathbb{D}}^\infty=\cap_{m\geq 1}{\mathbb{D}}^{n,\infty}$. Moreover, for $%
F\in ({\mathbb{D}}^{1,2})^{d},$ we let $\sigma _{F}$ denote the Malliavin
covariance matrix associated to $F:$%
\begin{equation*}
\sigma _{F}^{i,j}=\langle DF^{i},DF^{j}\rangle
=\sum_{k=1}^{N}\int_{0}^{1}D_{s}^{k}F^{i}D_{s}^{k}F^{j}ds,\quad i,j=1,\ldots
,d.
\end{equation*}%
If $\sigma _{F}$ is invertible, we denote through $\gamma _{F}$ the inverse
matrix. Finally, as usual, the notation $L$ will be used for the
Ornstein-Uhlenbeck operator and we recall that the Meyer inequality asserts
that $\Vert LF\Vert _{m,p}\leq C_{m,p}\Vert F\Vert _{m+2,p}$, for $F\in ({%
\mathbb{D}}^{m+2,\infty })^{d}$.

Our aim is to rephrase the results from the previous sections in the
framework of the Wiener space considered here. We introduce first the
localization random variable $\mathbf{\Theta .}$ We consider some random
variables $\Theta _{i}$ and some numbers $a_{i}>0,i=1,...,l+l^{\prime }$ and
we define
\begin{equation}
\mathbf{\Theta }=\prod_{i=1}^{l }\psi _{a_{i}}(\Theta _{i})\times
\prod_{i=l+1}^{l+l ^{\prime }}\phi _{a_{i}}(\Theta _{i}).  \label{W1}
\end{equation}%
with $\psi _{a_{i}},\phi _{a_{i}}$\ defined in (\ref{Mall10}). Following
what developed in Section \ref{th-abstract}, we define
\begin{equation*}
d{\mathbb{P}}_{\mathbf{\Theta }}=\mathbf{\Theta }d{\mathbb{P}}
\end{equation*}%
\ and
\begin{equation}
\Vert F\Vert _{p,\mathbf{\Theta }}^p={\mathbb{E}}_{\mathbf{\Theta }%
}(|F|^{p})\quad \mbox{and}\quad \Vert F\Vert _{l,p,\mathbf{\Theta }}^{p}={%
\mathbb{E}}_{\mathbf{\Theta }}(|F|_{l}^{p}).  \label{W2}
\end{equation}%
In the case $\mathbf{\Theta }=1$ we have $\Vert F\Vert _{p,\mathbf{\Theta }%
}=\Vert F\Vert _{p}$ and $\Vert F\Vert _{l,p,\mathbf{\Theta }}=\Vert F\Vert
_{l,p}.$ Moreover, given some $q\in {\mathbb{N}},p\geq 1,$ we denote
\begin{equation*}
\mathrm{m}_{q,p}(\mathbf{\Theta }):=1\vee \Vert \ln \mathbf{\Theta }\Vert
_{q,p,\mathbf{\Theta }}
\end{equation*}%
and%
\begin{equation}
\mathrm{U}_{q,p,\mathbf{\Theta }}(F):=\max \{1,{\mathbb{E}}_{\mathbf{\Theta }%
}((\det \sigma _{F})^{-p})(\left\Vert F\right\Vert _{q+2,p,\mathbf{\Theta }%
}+\left\Vert LF\right\Vert _{q,p,\mathbf{\Theta }})\}.  \label{W6}
\end{equation}%
In the case $\mathbf{\Theta }=1$ we have $\mathrm{m}_{q,p}(\mathbf{\Theta }%
)=1$ and
\begin{align}
\mathrm{U}_{q,p}(F) &:=\max \{1,{\mathbb{E}}((\det \sigma
_{F})^{-p})(\left\Vert F\right\Vert _{q+2,p}+\left\Vert LF\right\Vert
_{q,p})\}  \label{W7} \\
&\leq C\max \{1,{\mathbb{E}}((\det \sigma _{F})^{-p})\left\Vert F\right\Vert
_{q+2,p}\}  \notag
\end{align}%
the last inequality being a consequence of Meyer's inequality.

\smallskip

We rephrase now Theorem \ref{MainTh1}:

\begin{theorem}
\label{MainTh2} Let $q\in {\mathbb{N}}_{\ast }$. We consider the
localization random variable $\mathbf{\Theta }$ defined in (\ref{W1}) and we
assume that for every $p\in {\mathbb{N}}$ one has $\Vert \Theta _{i}\Vert
_{q+2,p}<\infty ,i=1,...,l+l^{\prime }.$ In particular $\mathrm{m}_{q+2,p}(%
\mathbf{\Theta })<\infty .$

\smallskip

\textbf{A}. Let $F\in ({\mathbb{D}}^{q+2,\infty })^d$ be such that $\mathrm{U%
}_{q,p,\mathbf{\Theta }}(F)<\infty .$ Then under ${\mathbb{P}}_{\mathbf{%
\Theta }}$ the law of $F$ is absolutely continuous with respect to the
Lebesgue measure. We denote by $p_{F,\mathbf{\Theta }}$ its density and we
have $p_{F,\mathbf{\Theta }}\in C^{q-1}({\mathbb{R}}^{d}).$ Moreover there
exist $C,a,b,p\in \mathcal{C}(q,d)$ such that for every $y\in {\mathbb{R}}%
^{d}$ and every multi index $\alpha =(\alpha _{1},...,\alpha _{k})\in
\{1,...,d\}^{k},k\in \{0,...,q\}$ one has%
\begin{equation}  \label{W8}
\left\vert \partial _{\alpha }p_{F,\mathbf{\Theta }}(y)\right\vert \leq C%
\mathrm{U}_{q,p,,\mathbf{\Theta }}^{a}(F)\times \mathrm{m}_{q+2,p}^{a}(%
\mathbf{\Theta })\times \big({\mathbb{P}}(\left\vert F-y\right\vert < 2 )%
\big)^{b}.
\end{equation}

\textbf{B}. Let $F,G\in ({\mathbb{D}}^{q+2,\infty })^d$ be such that $%
\mathrm{U}_{q+1,p,\mathbf{\Theta }}(F),\mathrm{U}_{q+1,p,\mathbf{\Theta }%
}(G)<\infty $ for every $p\in {\mathbb{N}}$ and let $p_{F,\mathbf{\Theta }}$
and $p_{G,\mathbf{\Theta }}$ be the densities of the laws of $F$
respectively of $G$ under ${\mathbb{P}}_{\mathbf{\Theta }}.$ Then there
exist $C,a,b,p\in \mathcal{C}(q,d)$ such that for every $y\in {\mathbb{R}}%
^{d}$ and every multi index $\alpha =(\alpha _{1},...,\alpha _{k})\in
\{1,2,...\}^{k},$ $0\leq k\leq q$ one has
\begin{equation}  \label{W9}
\begin{array}{rl}
\left\vert \partial _{\alpha }p_{F,\mathbf{\Theta }}(y)-\partial _{\alpha
}p_{G,\mathbf{\Theta }}(y)\right\vert \leq & C\mathrm{U}_{q+1,p,\mathbf{%
\Theta }}^{a}(F)\times \mathrm{U}_{q+1,p,\mathbf{\Theta }}^{a}(G)\times
\mathrm{m}_{q+2,p}^{a}(\mathbf{\Theta })\times \smallskip \\
& \times \big({\mathbb{P}}(\left\vert F-y\right\vert < 2 )+{\mathbb{P}}%
(\left\vert G-y\right\vert < 2 )\big)^{b}\times \smallskip \\
& \times \left( \Vert F-G\Vert _{q+2,p,\mathbf{\Theta }}+\Vert LF-LG\Vert
_{q,p,\mathbf{\Theta }}\right) .%
\end{array}%
\end{equation}
\end{theorem}

\begin{remark}
\label{use-mainW} The arguments used in Remark \ref{use-main} can be applied
here: the second factor in the estimates (\ref{W8}) and (\ref{W9}) can be
replaced, as $|y|>4$, with the queue of the law of $F$ and $G$. Also, by
using the Markov inequality, such factors can be over estimated by means of
any power of $(1+|y|)^{-1}$, for every $y\in{\mathbb{R}}^d$.
\end{remark}

\textbf{Proof}. One may prove Theorem \ref{MainTh2} just by repeating
exactly the same reasoning as in the proof of Theorem \ref{MainTh1}: all the
arguments are based on the properties of the norms from the finite
dimensional calculus and these properties are preserved in the infinite
dimensional case. However we give here a different proof: we obtain Theorem %
\ref{MainTh2} from Theorem \ref{MainTh1} by using a convergence argument.

We fix $n\geq 1$ and we denote%
\begin{equation*}
t_{n}^{k}=\frac{k}{2^{n}}\quad \mbox{and}\quad
I_{n}^{k}=(t_{n}^{k-1},t_{n}^{k}],\quad \mbox{with }k=1,\ldots ,2^{n}
\end{equation*}%
and%
\begin{eqnarray*}
\Delta _{n}^{k,i} &=&W^{i}(t_{n}^{k})-W^{i}(t_{n}^{k-1}),\qquad i=1,...,N, \\
\Delta _{n}^{k} &=&(\Delta _{n}^{k,1},...,\Delta _{n}^{k,N}),\qquad \Delta
_{n}=(\Delta _{n}^{1},...,\Delta _{n}^{2^{n}}).
\end{eqnarray*}%
We define the simple functionals of order $n$ to be the random variables
which are smooth functions of $\Delta _{n}$:
\begin{equation*}
\mathcal{S}_{n}=\{F=\phi (\Delta _{n}):\phi \in C_{p}^{\infty }({\mathbb{R}}%
^{N2^{n}})\}.
\end{equation*}%
And we define the simple processes of order $n$ by%
\begin{equation*}
\mathcal{P}_{n}=\Big\{U(s)=\sum_{i=1}^{N}%
\sum_{k=1}^{2^{n}}1_{I_{n}^{k}}(s)U_{k,i}\,:\, U_{k,i}\in \mathcal{S}_{n},\
k=1,\ldots ,2^{n}\Big\}.
\end{equation*}%
Our aim now is to identify the ``finite dimensional Malliavin calculus'' on $%
\mathcal{S}_{n}$ which gives us the standard Malliavin calculus on this
space. We take the basic random variable $V$ to be $\Delta _{n}=(\Delta
_{n}^{k,i},k=1,...,2^{n},i=1,...,n).$ The random variables $\Delta
_{n}^{k,i} $ are independent and
\begin{equation*}
p_{\Delta _{n}^{k,i}}(y)=\frac{1}{\sqrt{2\pi h_{n}}}e^{-y^{2}/2h_{n}},\qquad
h_{n}=\frac{1}{2^{n}}.
\end{equation*}%
In particular, with $y=(y_{k,i})_{i=1,...,N,k=1,...,2^{n}}$ we have $%
p_{\Delta _{n}}(y)=\prod_{i=1}^{N}\prod_{k=1}^{2^{n}}p_{\Delta
_{n}^{k,i}}(y_{k,i})$ and
\begin{equation*}
\partial _{y_{k,i}}\ln p_{\Delta _{n}}(y)=\frac{y_{k,i}}{h_{n}}.
\end{equation*}%
We take the weights $\pi _{k,i}=2^{-n/2}.$ Let $F\in \mathcal{S}_{n}\subset {%
\mathbb{D}}^{\infty }.$\ We recall that in the first section we have defined
the first order derivatives by
\begin{equation*}
D_{k,i}F=\pi _{k,i}\times \frac{\partial F}{\partial \Delta _{n}^{k,i}}
\end{equation*}%
and we used the norm%
\begin{align*}
\left\vert DF\right\vert _{1}^{2}
&=\sum_{i=1}^{N}\sum_{k=1}^{2^{n}}\left\vert D_{k,i}F\right\vert
^{2}=\sum_{i=1}^{N}\sum_{k=1}^{2^{n}}\pi _{k,i}^{2}\times \left\vert \frac{%
\partial F}{\partial \Delta _{n}^{k,i}}\right\vert ^{2} \\
&=\sum_{i=1}^{N}\sum_{k=1}^{2^{n}}\frac{1}{2^{n}}1_{I_{n}^{k}}(s)\left\vert
D_{s}^{i}F\right\vert ^{2}=\sum_{i=1}^{N}\int_{0}^{1}\left\vert
D_{s}^{i}F\right\vert ^{2}ds
\end{align*}%
where $D_{s}^{i}F$ is the standard Malliavin derivative. Similar
identification holds for the norms of higher order derivatives. Moreover we
recall that we defined%
\begin{eqnarray*}
LF &=&\delta (DF)=-\sum_{i=1}^{N}\sum_{k=1}^{2^{n}}(\partial _{\Delta
_{n}^{k,i}}(\pi _{k,i}D_{k,i}F)+\pi _{k,i}D_{k,i}F\partial _{\Delta
_{n}^{k,i}}(\ln p_{\Delta _{n}})(\Delta _{n}) \\
&=&-\sum_{i=1}^{N}\sum_{k=1}^{2^{n}}\frac{1}{2^{n}}\partial _{\Delta
_{n}^{k,i}}^{2}F+\sum_{i=1}^{N}\sum_{k=1}^{2^{n}}\frac{1}{2^{n}}\partial
_{\Delta _{n}^{k,i}}F\frac{\Delta _{k,i}}{2^{-n}} \\
&=&-\sum_{i=1}^{N}\sum_{k=1}^{2^{n}}\frac{1}{2^{n}}\partial _{\Delta
_{n}^{k,i}}^{2}F+\sum_{i=1}^{N}\sum_{k=1}^{2^{n}}\partial _{\Delta
_{n}^{k,i}}F\times \Delta _{k,i}
\end{eqnarray*}
and this is the Ornstein Uhlenbeck operator from the standard Malliavin
calculus. We conclude that the finite dimensional Malliavin calculus and the
standard Malliavin calculus coincide for simple functionals in $\mathcal{S}%
_{n}.$

We come now back to the proof of Theorem \ref{MainTh2} . We take $%
F_{n},G_{n},\Theta _{n,i}\in \mathcal{S}_{n},n\in {\mathbb{N}}%
,i=1,...,l+l^{\prime }$ which approximate $F,G,\Theta _{i}\in {\mathbb{D}}%
^{q+2,\infty },i=1,...,l+l^{\prime }.$We use Theorem \ref{MainTh1} for them
and then we pass to the limit in order to obtain the conclusion in Theorem %
\ref{MainTh2}. The fact that the constants which appear in Theorem \ref%
{MainTh1} belong to $\mathcal{C}(q,d)$, so do not depend on $n\in {\mathbb{N}%
},$ plays here a crucial role. $\square $

\begin{equation*}
\end{equation*}%
We give now a regularity property which is an easy consequence of the above
theorem.

\begin{theorem}
\textbf{A}. Let $F\in \D^{2,p},p>d$ such that $\P(\det \sigma _{F}>0)>0.$
Then, conditionally to $\{\det \sigma _{F}>0\},$ the law of $F$ is absolutely
continuous with respect to the Lebesgue measure and the density is lower
semi-continuous.

\smallskip

\textbf{B.} In particular the law of $F$ is locally lower bounded by the Lebesgue
measure $\lambda $ in the following sense: there exist an open set $%
D\subset \R^{d}$ and some $\delta >0$ such that for every Borel set $A$
\begin{equation*}
\P(F\in A)\geq \delta \lambda (A\cap D).
\end{equation*}%
\end{theorem}

\begin{remark}. The celebrated theorem of Bouleau and Hirsch \cite{bib:BH} says that if $%
F\in \D^{1,2}$ then, conditionally to $\{\det \sigma _{F}>0\},$ the law of $F$
is absolutely continuous. So it requires much less regularity than us. But
the new fact is that the conditional density is lower semi-continuous and in
particular is locally lower bounded by the Lebesgue measure. This last
property turns out to be especially interesting - see the joint paper \cite{bib:CLT}.
\end{remark}

\textbf{Proof}. For $\varepsilon >0$ we consider the localization function $%
\psi _{\varepsilon }$ defined in (\ref{Mall10}) and we denote $\Theta
_{\varepsilon }=\psi _{\varepsilon }(\det \sigma _{F}).$ By Theorem \ref{MainTh2} we know that under $\P_{\Theta _{\varepsilon }}$ the law of $F$ is
absolutely continuous and has a continuous density $p_{\Theta _{\varepsilon
}}.$ Let $A$ be a Borel set with $\lambda (A)=0$ where $\lambda $ is the
Lebesgue measure. Since $\Theta _{\varepsilon }\uparrow \Theta :=1_{\{\det
\sigma _{F}>0\}}$ we have%
\begin{align*}
\P_{\Theta}(F\in A)
&=\P_{\{\det \sigma _{F}>0\}}(F \in A)
=\frac{1}{\P(\det \sigma _{F}>0)}\E(1_{\{F\in A\}}\times \Theta ) \\
&=\frac{1}{\P(\det \sigma _{F}>0)}\lim_{\varepsilon \rightarrow
0}\E(1_{\{F\in A\}}\times \Theta _{\varepsilon })=0.
\end{align*}%
So we may find $p_{\Theta }$ such that
\begin{equation*}
\E(f(F)\Theta )=\int f(x)p_{\Theta }(x)dx.
\end{equation*}%
For $f\geq 0$ we have
\begin{equation*}
\int f(x)p_{\Theta _{\varepsilon }}(x)dx=\E(f(F)\Theta _{\varepsilon })\leq
\E(f(F)\Theta )=\int f(x)p_{\Theta }(x)dx,
\end{equation*}%
so that $p_{\Theta }\geq p_{\Theta _{\varepsilon }}$ a.e. This implies that $p_{\Theta }\geq \sup_{\varepsilon>0}p_{\Theta _{\varepsilon }}$. We claim that
\begin{equation*}
p_{\Theta }=\sup_{\varepsilon >0}p_{\Theta _{\varepsilon }}
\end{equation*}%
which gives that $p_{\Theta }$ is lower semi-continuous. In fact, set $A=\{x: p_{\Theta }(x)>\sup_{\varepsilon >0}p_{\Theta _{\varepsilon
}}(x)\}$. If $%
\lambda (A)>0$ then we may find $\delta >0$ such that $\lambda (A_{\delta })>0$
with $A_{\delta }=\{x:p_{\Theta }(x)>\delta +\sup_{\varepsilon >0}p_{\Theta
_{\varepsilon }}(x)\}.$ Then%
\begin{align*}
\int_{A_{\delta }}p_{\Theta }(x)dx
&=\P_{\{\det \sigma _{F}>0\}}(F\in A_{\delta })
=\frac{1}{\P(\det \sigma _{F}>0)}\lim_{\varepsilon \to
0}\E(1_{\{F\in A_{\delta }\}}\times \Theta _{\varepsilon }) \\
&=\lim_{\varepsilon \rightarrow 0}\int_{A_{\delta }}p_{\Theta _{\varepsilon
}}(x)dx\leq \int_{A_{\delta }}(p_{\Theta }(x)-\delta )dx
\end{align*}%
and this would give $\lambda (A_{\delta })=0.$

\smallskip

The assertion \textbf{B} is immediate: since $p_{\Theta }=\sup_{\varepsilon
>0}p_{\Theta _{\varepsilon }}$ is not identically null we may find $%
\varepsilon >0$ and $x_{0}\in \R^{d}$ such that $p_{\Theta _{\varepsilon
}}(x_{0})>0.$ And since $p_{\Theta _{\varepsilon }}$ is a continuous
function we may find $r,\delta >0$ such that $p_{\Theta _{\varepsilon
}}(x)\geq \delta $ for $x\in B_{r}(x_{0}).$ It follows that
\begin{align*}
\P(F \in A)
&\geq \P(\{F\in A\}\cap \{F\in B_{r}(x_{0})\}\cap \{\sigma
_{F}>0\}) \\
&=\P(\sigma _{F}>0)\int_{A\cap B_{r}(x_{0})}p_{\Theta }(x)dx \\
&\geq \P(\sigma _{F}>0)\int_{A\cap B_{r}(x_{0})}p_{\Theta _{\varepsilon
}}(x)dx\geq \delta \P(\sigma _{F}>0)\lambda (A\cap B_{r}(x_{0})).
\end{align*}%
$\square $

\medskip

We rephrase now other consequences of Theorem \ref{MainTh2}. We begin with
the regularization Lemma \ref{Regularization}. We recall that $\gamma
_{\delta }$ is the centred Gaussian density with variance $\delta >0.$

\begin{lemma}
\label{RegularizationW}There exist some universal constants $C,p,a\in
\mathcal{C}(d)$ such that for every $\varepsilon >0,\delta >0$ and every $%
F\in ({\mathbb{D}}^{3,\infty })^{d}$\ one has%
\begin{equation}  \label{W10}
\left\vert {\mathbb{E}}(f(F))-{\mathbb{E}}(f\ast \gamma _{\delta
}(F))\right\vert \leq C\left\Vert f\right\Vert _{\infty }\big({\mathbb{P}}%
(\sigma _{F}<\varepsilon )+\frac{\sqrt{\delta}}{\varepsilon ^{p}}(1+ \Vert
F\Vert _{3,p})^{a}\big)
\end{equation}%
for every $f\in C_{b}({\mathbb{R}}^{d}).$
\end{lemma}

\textbf{Proof.} The proof is identical with the one of Lemma \ref%
{Regularization} so we skip it (an approximation procedure may also been
used). We mention that due to Meyer's inequalities $\left\Vert LF\right\Vert
_{1,p}$ does no more appear here. $\square $

\medskip

We consider now the distances $d_{m}$ defined in (\ref{Main9}) and we
rewrite Theorem \ref{Distance}:

\begin{theorem}
\label{Distance copy(1)} Let $k\in {\mathbb{N}}.$ There exist universal
constants $C,p,b\in \mathcal{C}(d,k)$ such that for every $F,G\in ({\mathbb{D%
}}^{3,\infty })^{d}$ and every $\varepsilon >0$ one has%
\begin{equation}  \label{W11}
d_{0}(F,G) \leq \displaystyle\frac{C}{\varepsilon ^{b}}(1+\|F\|_{3,p}+\|G%
\|_{3,p}\|)^bd_{k}^{\frac{1}{k+1}}(F,G) +C{\mathbb{P}}(\det \sigma _{F}
<\varepsilon )+C{\mathbb{P}}(\det \sigma _{G}<\varepsilon ).
\end{equation}
\end{theorem}

\textbf{Proof.} The proof is identical with the one of Theorem \ref{Distance}
so we skip it. $\square $

\medskip

We give now the convergence results.

\begin{theorem}
\label{GeneralW}We consider a sequence of functionals $%
F_{n}=(F_{n,1},...,F_{n,d})\in ({\mathbb{D}}^{3,\infty })^{d},n\in {\mathbb{N%
}}$ and we assume that%
\begin{equation}
\begin{array}{rcl}
i) &  & \quad \sup_{n}\left\Vert F_{n}\right\Vert _{3,p}<\infty ,\qquad
\forall p\geq 1,\smallskip  \\
ii) &  & \quad \displaystyle\limsup_{\varepsilon \rightarrow
0}\limsup_{n\rightarrow \infty }{\mathbb{P}}(\det \sigma
_{F_{n}}<\varepsilon )=0.%
\end{array}
\label{W12}
\end{equation}%
Suppose also that $\lim_{n}F_{n}=F$ in distribution and $%
\lim_{n}\E(F_{n})=\E(F)$. Then
\begin{equation*}
\lim_{n}d_{TV}(F,F_{n})=0.
\end{equation*}%
In particular if the laws of $F$ and $F_{n}$ are absolutely continuous with
densities $p_{F}$ and $p_{F_{n}}$ then
\begin{equation*}
\lim_{n}\int \left\vert p_{F}(x)-p_{F_{n}}(x)\right\vert dx=0.
\end{equation*}
\end{theorem}

\textbf{Proof.} The proof is identical with the one of Theorem \ref{General}
so we skip it. $\square $

\medskip

In the framework of Wiener functionals we are able to obtain one more result:

\begin{corollary}
\label{corW} Let $F_{n}\in ({\mathbb{D}}^{3,\infty })^{d},n\in {\mathbb{N}}$
such that $\sup_{n}\left\Vert F_{n}\right\Vert _{3,p}<\infty $ for every $%
p\in {\mathbb{N}}$. Consider also $F\in ({\mathbb{D}}^{2,2})^{d}$ such that $%
\det \sigma _{F}>0$ almost surely. If $\lim_{n}F_{n}=F$ in $L^{2}$ then $%
\lim_{n}F_{n}=F$ in $d_{TV}$.
\end{corollary}

\textbf{Proof}. We will prove that $\lim_{n}\left\langle
DF_{n}^{i},DF_{n}^{j}\right\rangle =\left\langle DF^{i},DF^{j}\right\rangle $
in probability for every $i,j=1,...,d.$ These implies $\lim_{n}\det \sigma
_{F_{n}}=\det \sigma _{F}$ in probability so that $\limsup_n {\mathbb{P}}%
(\det\sigma_{F_n}<\varepsilon)\leq P(\det \sigma _{F}<2\varepsilon ).$ And
since $\lim_{\varepsilon \to 0}P(\det \sigma _{F}>\varepsilon )=1$, we
obtain (\ref{W12}) $ii)$ and this permits to conclude by applying Theorem %
\ref{GeneralW}.

We denote by $f_{k},k\in {\mathbb{N}}$, respectively by $f_{k,n},k\in {%
\mathbb{N}}$, the kernels of the chaos expansion of $F$, respectively of $%
F_{n}.$ So we have%
\begin{equation*}
F=\sum_{k=0}^{\infty }I_{k}(f_{k})\quad \mbox{and}\quad
F_{n}=\sum_{k=0}^{\infty }I_{k}(f_{k,n})
\end{equation*}%
where $I_{k}$ denotes the multiple integral of order $k.$ For $N\in {\mathbb{%
N}}$ we write $F=S_{N}+R_{N}$ with $S_{N}=\sum_{k=0}^{N}I_{k}(f_{k})$ and $%
R_{N}=\sum_{k=N+1}^{N}I_{k}(f_{k})$. With similar notations, we set $%
F_{n}=S_{N,n}+R_{N,n}.$ We write%
\begin{equation*}
\left\vert \left\langle DF_{n}^{i},DF_{n}^{j}\right\rangle -\left\langle
DF^{i},DF^{j}\right\rangle \right\vert \leq (\left\vert DF\right\vert
+\left\vert DF_{n}\right\vert )\times \left\vert D(F-F_{n})\right\vert .
\end{equation*}%
Since the sequence $\left\vert DF_{n}\right\vert ,n\in {\mathbb{N}}$ is
bounded in $L^{1}$ our conclusion follows as soon as we check that $%
\lim_{n}\left\vert D(F-F_{n})\right\vert =0$ in probability. We fix $%
\varepsilon >0$ and we write%
\begin{align*}
{\mathbb{P}}(\left\vert D(F-F_{n})\right\vert \geq \varepsilon ) &\leq {%
\mathbb{P}}\Big(\left\vert D(S_{N}-S_{N,n})\right\vert \geq \frac{1}{2}%
\varepsilon \Big)+{\mathbb{P}}\Big(\left\vert D(R_{N}-R_{N,n})\right\vert
\geq \frac{1}{2}\varepsilon \Big) \\
&=:I_{N,\varepsilon ,n}+J_{N,\varepsilon ,n}.
\end{align*}%
Using Chebyshev's inequality
\begin{equation*}
J_{N,\varepsilon ,n}\leq \frac{4}{\varepsilon ^{2}}{\mathbb{E}}\left\vert
D(R_{N}-R_{N,n})\right\vert ^{2}\leq \frac{8}{\varepsilon ^{2}}({\mathbb{E}}%
\left\vert DR_{N}\right\vert ^{2}+{\mathbb{E}}\left\vert DR_{N,n}\right\vert
^{2}).
\end{equation*}%
Since ${\mathbb{E}}\left\vert DI_{k}(f_{k})\right\vert ^{2}=k\times
k!\left\Vert f_{k}\right\Vert _{L^{2}[0,T]^{k}}^{2}$ we obtain%
\begin{align*}
{\mathbb{E}}\left\vert DR_{N}\right\vert ^{2} &=\sum_{k=N+1}^{\infty
}k\times k!\left\Vert f_{k}\right\Vert _{L^{2}[0,1]^{k}}^{2}\leq \frac{1}{N+2%
}\sum_{k=N+1}^{\infty }(k+1)\times k\times k!\left\Vert f_{k}\right\Vert
_{L^{2}[0,1]^{k}}^{2} \\
&\leq \frac{1}{N+2}{\mathbb{E}}(\vert D^{2}F\vert ^{2})\leq \frac{1}{N+2}%
\left\Vert F\right\Vert _{2,2}^2
\end{align*}%
and a similar inequality holds for ${\mathbb{E}}\left\vert
DR_{N,n}\right\vert ^{2}.$ We conclude that
\begin{equation*}
J_{N,\varepsilon ,n}\leq \frac{8}{(N+2)\varepsilon ^{2}}(\left\Vert
F\right\Vert _{2,2}+\sup_{n}\left\Vert F_{n}\right\Vert _{2,2}).
\end{equation*}%
Moreover, since $\lim_{n}\left\Vert F-F_{n}\right\Vert _{2}=0$ we have $%
\lim_{n}\left\Vert f_{k}-f_{k,n}\right\Vert _{L^{2}[0,T]^{k}}^{2}=0$ for
every $k\in {\mathbb{N}}$ and this implies $\lim_{n}{\mathbb{E}}\left\vert
D(S_{N}-S_{N,n})\right\vert ^{2}=0.$ Finally $\limsup _{n}I_{N,\varepsilon
,n}=0$ for each fixed $N$ and $\varepsilon .$ So we obtain
\begin{equation*}
\limsup_n{\mathbb{P}}(\left\vert D(F-F_{n})\right\vert \geq \varepsilon
)\leq \frac{8}{N^{2}\varepsilon ^{2}}(\left\Vert F\right\Vert
_{2,2}+\sup_{n}\left\Vert F_{n}\right\Vert _{2,2}).
\end{equation*}%
Since this is true for each $N$ the above limit is null. $\square $

\begin{remark}
As an immediate consequence of Corollary \ref{corW} one may obtain the
following result. Let $X_{t}$ be a diffusion process with coefficients in $%
C_{b}^{\infty }$ and suppose that the weak H\"{o}rmander condition holds in $%
x=X_{0}.$ Consider also the Euler scheme \ $X_{t}^{n}$ of step $\frac{1}{n}.$
Then for every $q\in {\mathbb{N}}$ one has $d_{-q}(\mu _{X_{t}},\mu
_{X_{t}^{n}})\rightarrow 0.$ This type of result has already been obtained
in \cite{[BT1]}, \cite{bib:[BT2]} and in \cite{[Gu]}: there, under more
restrictive assumptions (uniform H\"{o}rmander condition) one obtain the
above result and moreover, one gives a development in Taylor series of the
error.
\end{remark}

\section{Proof of Theorem \protect\ref{MainTh1}}

\label{th-proof}

This section is devoted to the proof of Theorem \ref{MainTh1}. We are in the
framework defined in Section \ref{th-abstract} and we use all the notation
introduced there. In the following subsection we recall and develop some
basic results concerning integration by parts formulas from \cite{bib:[BCl]}.

\subsection{Integration by parts formulae}

\label{sect-dualJ}

By using standard integration by parts formulas, one gets the duality
between $\delta $ and $D$ and the standard computation rules (see \cite%
{bib:[BCl]}, Proposition 1 and Lemma 1): under our assumption, for every $%
F\in\mathcal{S}^d$, $U\in \mathcal{P}$ and $\phi :\mathbb{R}^{d}\rightarrow
\mathbb{R}$ smooth,
\begin{align}
{\mathbb{E}}(\left\langle DF,U\right\rangle _{J})&={\mathbb{E}}(F\delta (U)),
\label{1.2} \\
D\phi (F)&=\sum_{r=1}^{d}\partial _{r}\phi (F)DF^{r},  \label{1.3} \\
\delta (FU)&=F\delta (U)-\left\langle DF,U\right\rangle_J,  \label{1.4} \\
L \phi(F)&=\sum_{r=1}^{d}\partial _{r}\phi (F)LF^{r}
-\sum_{r,r^{\prime}=1}^{d}\partial _{r,r^{\prime}}\phi (F)\langle DF^{r},
DF^{r^{\prime}}\rangle_J,  \label{CRL}
\end{align}
with the understanding $d=1$ in \eqref{1.2} and \eqref{1.4}. Once the above
equalities are done, the integration by parts formulas can be stated (see
\cite{bib:[BCl]}, Theorem 1 and 2):

\begin{theorem}
\label{IPP} Let $F\in \mathcal{S}^{d}$ be such that
\begin{equation}
{\mathbb{E}}(\left\vert \det \gamma _F\right\vert ^{p})<\infty \quad \forall
p\geq 1,  \label{1.5}
\end{equation}%
$\gamma_F$ denoting the inverse of the Malliavin covariance matrix $\sigma_F$%
. Then, for every $G\in\mathcal{S}$ and for every smooth function $\phi\in
C^\infty_b({\mathbb{R}}^d)$ one has
\begin{equation}
{\mathbb{E}}\left( \partial _{r}\phi (F)G\right) ={\mathbb{E}}\left( \phi
(F)H_{r}(F,G)\right),\quad r=1,...,d,  \label{IPP1}
\end{equation}%
with
\begin{equation}  \label{weight1}
H_{r}(F,G) =\displaystyle\sum_{r^{\prime }=1}^{d}\delta
(G\gamma_F^{r^{\prime },r}DF^{r^{\prime }}) =\displaystyle\sum_{r^{\prime
}=1}^{d}\big( G\delta (\gamma_F ^{r^{\prime },r}DF^{r^{\prime }})-\gamma_F
^{r^{\prime },r}\langle DF^{r^{\prime }},DG\rangle _{J}\big) .
\end{equation}
Moreover, for every $q\in{\mathbb{N}}^*$ and multi-index $\beta =(\beta
_{1},\ldots ,\beta _{q})\in \{1,\ldots ,d\}^{q}$ then
\begin{equation}
{\mathbb{E}}\left( \partial _{\beta }\phi (F)G\right) ={\mathbb{E}}\left(
\phi (F)H_{\beta }^{q}(F,G)\right)  \label{IPP2}
\end{equation}%
where the weights $H^{q}_\beta(F,G)$ are defined recursively by (\ref%
{weight1}) if $q=1$ and for $q>1$,
\begin{equation}
H_{\beta }^{q}(F,G)=H_{\beta _{1}}\big( F,H_{(\beta _{2},\ldots ,\beta
_{q})}^{q-1}(F,G)\big) .  \label{weight2}
\end{equation}
\end{theorem}

\subsection{Estimates of the weights}

In this section we give estimates of the weights $H_{\alpha }^{q}(F,G)$
appearing in the integration by parts formulae of Theorem \ref{IPP} using
the norms introduced in (\ref{NORM1}).\ We first deal with useful estimates
for the inverse of the Malliavin covariance matrix. For $F\in \mathcal{S}%
^{d} $, we set
\begin{equation}
\m F=\max \Big(1,\frac{1}{\det \sigma _{F}}\Big).  \label{FM}
\end{equation}

\begin{proposition}
\label{Bgamma} \textbf{A}. If $F\in \mathcal{S}^{d}$ then $\forall l\in
\mathbb{N}$ one has

\begin{equation}
|\gamma _{F}|_{l}\leq C_{l,d}{\m F}^{l+1}(1+|F|_{1,l+1}^{2d(l+1)}).
\label{gamma1}
\end{equation}%
\textbf{B.} If $F,\overline{F}\in \mathcal{S}^{d}$ then $\forall l\in
\mathbb{N}$ one has
\begin{equation}
|\gamma _{F}-\gamma _{\overline{F}}|_{l}\leq C_{l,d}\m F^{l+1}\,%
\m{\overline{F}}^{l+1}(1+|F|_{1,l+1}+|\overline{F}|_{1,l+1})^{2d(l+3)}\left%
\vert F-\overline{F}\right\vert _{1,l+1}.  \label{gamma2}
\end{equation}
\end{proposition}

\textbf{Proof.} \textbf{A} is proved in \cite{bib:[BCl]}, Proposition 2. As
for \textbf{B}, we use the following estimates proved in \cite{bib:[BCl]}
(see Lemma 2 and the proof of Proposition 2):
\begin{align}
|\left\langle DF,DG\right\rangle _{J}|_{l}& \leq 2^{l}\sum_{l_{1}+l_{2}\leq
l}|F|_{1,l_{1}+1}|G|_{1,l_{2}+1},  \label{scal1} \\
|F\times G|_{l}& \leq 2^{l}\sum_{l_{1}+l_{2}\leq l}|F|_{l_{1}}|G|_{l_{2}},
\label{prod} \\
|(\det \sigma _{F})^{-1}|_{l}& \leq C_{l_{1}}\m %
F^{l+1}(1+|F|_{1,l+1}^{2l_{1}d}),  \label{est-det} \\
|\det \sigma _{F}|_{l}& \leq C_{l}(1+|F|_{1,l+1}^{2l_{1}d})  \label{det}
\end{align}%
So, by (\ref{scal1}), we have
\begin{equation*}
\left\vert \sigma _{F}^{r,r^{\prime }}-\sigma _{\overline{F}}^{r,r^{\prime
}}\right\vert _{l}\leq C_{l,d}\left\vert F-\overline{F}\right\vert
_{1,l+1}(\left\vert F\right\vert _{1,l+1}+\left\vert \overline{F}\right\vert
_{1,l+1})
\end{equation*}%
and then, by (\ref{prod}) and  (\ref{det})%
\begin{align}
\left\vert \det \sigma _{F}-\det \sigma _{\overline{F}}\right\vert _{l}&
\leq C_{l,d}\left\vert F-\overline{F}\right\vert _{1,l+1}(\left\vert
F\right\vert _{1,l+1}+\left\vert \overline{F}\right\vert _{1,l+1})^{2d-1}
\label{b} \\
\left\vert \widehat{\sigma }_{F}^{r,r^{\prime }}-\widehat{\sigma }_{%
\overline{F}}^{r,r^{\prime }}\right\vert _{l}& \leq C_{l,d}\left\vert F-%
\overline{F}\right\vert _{1,l+1}(\left\vert F\right\vert _{1,l+1}+\left\vert
\overline{F}\right\vert _{1,l+1})^{2d-3},  \notag
\end{align}%
in which $\widehat{\sigma }_{F}$ denotes the algebraic complement of $\sigma
_{F}$. Then, by using also \eqref{est-det}
\begin{align*}
\left\vert (\det \sigma _{F})^{-1}-(\det \sigma _{\overline{F}%
})^{-1}\right\vert _{l}\leq & C_{l,d}\left\vert (\det \sigma
_{F})^{-1}\right\vert _{l}\left\vert (\det \sigma _{\overline{F}%
})^{-1}\right\vert _{l}\times  \\
& \times \left\vert \det \sigma _{F}-\det \sigma _{\overline{F}}\right\vert
_{l} \\
\leq & C_{l,d}\m F^{l+1}\m{\overline{F}}^{l+1}|F-\overline{F}|_{1,l+1}\big(%
1+|F|_{1,l}+|\overline{F}|_{1,l}\big)^{2(l+1)d}.
\end{align*}%
Since $\gamma ^{r,r^{\prime }}(F)=(\det \sigma _{F})^{-1}\widehat{\sigma }%
^{r,r^{\prime }}(F)$, \eqref{gamma2} follows by using the above estimates. $%
\square $\medskip

We define now%
\begin{equation}
L_{r}^{\gamma }(F)=\sum_{r^{\prime }=1}^{d}\delta (\gamma _{F}^{r^{\prime
},r}DF^{r^{\prime }})=\sum_{r^{\prime }=1}^{d}\big(\gamma _{F}^{r^{\prime
},r}LF^{r^{\prime }}-\<D\gamma _{F}^{r^{\prime },r},DF^{r^{\prime }}\>\big).
\label{gamma3}
\end{equation}%
Using (\ref{gamma2}) one can easily check that for every $l\in {\mathbb{N}}$%
\begin{equation}
\left\vert L_{r}^{\gamma }(F)\right\vert _{l}\leq C_{l,d}\m F^{l+2}\big(%
1+\left\vert F\right\vert _{1,l+1}^{2d(l+2)}+\left\vert L(F)\right\vert
_{l}^{2}\big).  \label{gamma4}
\end{equation}%
And by using both (\ref{gamma2}) and \eqref{gamma3}, one immediately gets
\begin{equation}
|L_{r}^{\gamma }(F)-L_{r}^{\gamma }(\overline{F})|_{l}\leq C_{l,d}Q_{l}(F,%
\overline{F})(|F-\overline{F}|_{1,l+2}+|L(F-\overline{F}))|_{l}),
\label{gamma5}
\end{equation}%
where
\begin{equation}
Q_{l}(F,\overline{F})=\m F^{l+2}\m{\overline{F}}^{l+2}\big(1+\left\vert
F\right\vert _{1,l+2}^{2d(l+4)}+\left\vert L(F)\right\vert
_{l}^{2}+\left\vert \overline{F}\right\vert _{1,l+2}^{2d(l+4)}+\left\vert L(%
\overline{F})\right\vert _{l}^{2}\big).  \label{gamma6}
\end{equation}%
For $F\in \mathcal{S}^{d}$, we define the linear operator $T_{r}(F,\cdot ):%
\mathcal{S}\rightarrow \mathcal{S},r=1,...,d$ by
\begin{equation*}
T_{r}(F,G)=\left\langle DG,(\gamma _{F}DF)^{r}\right\rangle ,
\end{equation*}%
where $(\gamma _{F}DF)^{r}=\sum_{r^{\prime }=1}^{d}\gamma _{F}^{r^{\prime
},r}DF^{r^{\prime }}$. Moreover, for a multi-index $\beta =(\beta
_{1},..,\beta _{q})$ we denote $\left\vert \beta \right\vert =q$ and we
define by induction
\begin{equation*}
T_{\beta }(F,G)=T_{\beta _{q}}(F,T_{(\beta _{1},...,\beta _{q-1})}(F,G)).
\end{equation*}

For $l\in {\mathbb{N}}$ and $F,\overline{F}\in \mathcal{S}^{d}$, we denote%
\begin{align}
& \Theta _{l}(F)=\m F^{l}(1+\left\vert F\right\vert
_{1,l+1}^{2d(l+1)})\qquad \mbox{and}  \label{gamma7} \\
& \Theta _{l}(F,\overline{F})=\m F^{l}\m{\overline{F}}^{l}\big(%
1+|F|_{1,l}^{2d(l+2)}+|\overline{F}|_{1,l}^{2d(l+2)}\big).  \label{gamma8}
\end{align}%
We notice that $\Theta _{l}(F)\leq \Theta _{l}(F,\overline{F})$, $\Theta
_{l}(F)\leq \Theta _{l+1}(F)$ and $\Theta _{l}(F,\overline{F})\leq \Theta
_{l+1}(F,\overline{F})$.

\begin{proposition}
Let $F,\overline{F}\in \mathcal{S}^{d}$ and $G,\overline{G}\in \mathcal{S}.$
Then for every $l\in {\mathbb{N}}$ and for every multi index $\beta $ with $%
|\beta |=q\geq 1$ one has
\begin{equation}
\left\vert T_{\beta }(F,G)\right\vert _{l}\leq C\Theta
_{l+q}^{q}(F)\,|G|_{1,l+q}  \label{gamma9}
\end{equation}%
and%
\begin{eqnarray}
\left\vert T_{\beta }(F,G)-T_{\beta }(\overline{F},\overline{G})\right\vert
_{l} &\leq &C\Theta _{l+q}^{\frac{q(q+1)}{2}}(F,\overline{F})\big(%
1+|G|_{1,l+q}+|\overline{G}|_{1,l+q}\big)^{q}\times   \label{gamma10} \\
&&\times (\left\vert F-\overline{F}\right\vert _{1,l+q}+\left\vert G-%
\overline{G}\right\vert _{1,l+q})  \notag
\end{eqnarray}%
where $C\in \mathcal{C}(l,d,q)$.
\end{proposition}

\textbf{Proof.} \eqref{gamma9} follows from \cite{bib:[BCl]} (see (26) in
the proof of Theorem 3 therein). We prove \eqref{gamma10} by recurrence.
Hereafter, $C$ denotes a constant in $\mathcal{C}(l,d,q)$, possibly varying
from a line to another. For $|\beta |=q=1$ we have
\begin{eqnarray*}
\left\vert T_{r}(F,G)-T_{r}(\overline{F},\overline{G})\right\vert _{l} &\leq
&\left\vert G-\overline{G}\right\vert _{l+1}(\left\vert F\right\vert
_{1,l+1}+\left\vert \overline{F}\right\vert _{1,l+1})(\left\vert \gamma
_{F}\right\vert _{l}+\left\vert \gamma _{\overline{F}}\right\vert _{l})+ \\
&&+\left\vert F-\overline{F}\right\vert _{1,l+1}(\left\vert G\right\vert
_{1,l+1}+\left\vert \overline{G}\right\vert _{1,l+1})(\left\vert \gamma
_{F}\right\vert _{l}+\left\vert \gamma _{\overline{F}}\right\vert _{l})+ \\
&&+\left\vert \gamma _{F}-\gamma _{\overline{F}}\right\vert _{l}(\left\vert
F\right\vert _{1,l+1}+\left\vert \overline{F}\right\vert
_{1,l+1})(\left\vert G\right\vert _{1,l+1}+\left\vert \overline{G}%
\right\vert _{1,l+1}).
\end{eqnarray*}%
Using \eqref{gamma2}, we obtain (\ref{gamma10}). The case $|\beta |=q\geq 2$%
, easily follows by induction. $\square $\medskip

We can now establish estimates for the weights $H^{q}$. For $l\geq 1$ and $F,%
\overline{F}\in \mathcal{S}^{d}$, we set (with the understanding $|\cdot
|_{0}=|\cdot |$)
\begin{align}
& A_{l}(F)=\m F^{l+1}\big(1+|F|_{1,l+1}^{2d(l+2)}+|LF|_{l-1}^{2}\big)
\label{gamma11'} \\
& A_{l}(F,\overline{F})=\m F^{l+1}\m {\overline{F}}^{l+1}\big(%
1+|F|_{1,l+1}^{2d(l+3)}+|\overline{F}|_{1,l+1}^{2d(l+3)}+|LF|_{l-1}^{2}+|L%
\overline{F}|_{l-1}^{2}\big).  \label{AAn}
\end{align}

\begin{theorem}
\label{thHq} \textbf{A}. For $l\in {\mathbb{N}}$ and $q\in \mathbb{N}^{\ast}
$ there exists $C\in \mathcal{C}(l,d,q)$ such that for every $F\in \mathcal{S%
}^{d}$, $G\in \mathcal{S}$ and for every multi-index $\beta =(\beta
_{1},..,\beta _{q})$
\begin{equation}
\left\vert H_{\beta }^{q}(F,G)\right\vert _{l}\leq CA_{l+q}(F)^{q}|G|_{l+q},
\label{gamma11}
\end{equation}%
$A_{l+q}(F)$ being defined in \eqref{gamma11'}.

\smallskip

\textbf{B}. There exists $C\in \mathcal{C}(l,d,q)$ such that for every $F,%
\overline{F}\in \mathcal{S}^{d}$ , $G,\overline{G}\in \mathcal{S}$ and every
multi-index $\beta =(\beta _{1},..,\beta _{q})$
\begin{equation}
\begin{array}{rl}
\left\vert H_{\beta }^{q}(F,G)-H_{\beta }^{q}(\overline{F},\overline{G}%
)\right\vert _{l}\leq & CA_{l+q}(F,\overline{F})^{\frac{q(q+1)}{2}}\times %
\big(1+|G|_{l+q}+|\overline{G}|_{l+q}\big)^{q}\times \smallskip \\
& \times \big(\left\vert F-\overline{F}\right\vert _{l+q+1}+\left\vert L(F-%
\overline{F})\right\vert _{l+q-1}+\left\vert G-\overline{G}\right\vert _{l+q}%
\big).%
\end{array}
\label{gamma12}
\end{equation}
\end{theorem}

\textbf{Proof.} \textbf{A.} Suppose $q=1$ and $\beta =r$. Then,
\begin{equation*}
|H_{r}(F,G)|_{l}\leq C|G|_{l}|L_{r}^{\gamma }(F)|_{l}+C|T_{r}(F,G)|_{l}.
\end{equation*}%
By using \eqref{gamma4} and \eqref{gamma9}, we can write
\begin{align*}
|H_{r}(F,G)|_{l}\leq & C\m F^{l+2}\big(1+|F|_{1,l+1}^{2d(l+2)}+|LF|_{l}^{2}%
\big)|G|_{l}+\m F^{l+1}\big(1+|F|_{1,l+2}^{2d(l+3)}\big)|G|_{l+1} \\
\leq & A_{l+1}(F)|G|_{l+1}.
\end{align*}%
So, the statement holds for $q=1$. And for $q>1$ it follows by iteration and
by using the fact that $A_{l+1}(F)\leq A_{l+q}(F)$.

\smallskip

\textbf{B.} Suppose $q=1$ and $\beta=r$. Then,
\begin{align*}
|H_r(F,G)-H_r(\overline{F}, \overline{G})|_l \leq & C\big( %
|L_r^{\gamma}(F)|_l|G-\overline{G}|_l +|\overline{G}|_l|L_r^{%
\gamma}(F)-L_r^{\gamma}(\overline{F})|_l + \\
&\quad+|T_r(F,G)-T_r(\overline{F},\overline{G})|_l\big).
\end{align*}
Now, estimate \eqref{gamma12} follows by using \eqref{gamma4} and %
\eqref{gamma5}, \eqref{gamma10}. In the iteration for $q>1$, it suffices to
observe that $A_l(F)\leq A_l(F,\overline{F})$, $A_l(F)\leq A_{l+1}(F)$ and $%
A_l(F,\overline{F})\leq A_{l+1}(F,\overline{F})$. $\square$\medskip

\subsection{Localized representation formulas for the density}

\label{sect-locIBP}

In this section we discuss localized integration by parts formulas with the
localization random variable $\mathbf{\Theta }$ defined in (\ref{Mall10'}).
We will use the norms $\Vert F\Vert _{p,\mathbf{\Theta }}$ and $\Vert F\Vert
_{1,l,p,\mathbf{\Theta }},\Vert F\Vert _{l,p,\mathbf{\Theta }}$ defined in (%
\ref{NORM1}). We also recall that $\mathrm{m}_{q,p}(\mathbf{\Theta })$ is
defined in (\ref{Mall1}) and that an estimate of this quantity is given in (%
\ref{Mall11}).

We give now the integration by parts formula with respect to ${\mathbb{P}}_{%
\mathbf{\Theta }}$ (that is, locally) and we study the regularity of the law
starting from the results in \cite{bib:[BC0]}.

Once for all, in addition to $\mathrm{m}_{q,p}(\mathbf{\Theta })$ we define
the following quantities: for $p\geq 1$, $q\in {\mathbb{N}}$, $F\in \mathcal{%
S}^{d}$,%
\begin{equation}
\begin{array}{rcl}
{\mathrm{S}}_{F,\mathbf{\Theta }}(p) & = & \max \{1,\left\Vert (\det \sigma
_{F})^{-1}\right\Vert _{p,\mathbf{\Theta }}\}, \\
{\mathrm{Q}}_{F,\mathbf{\Theta }}(q,p) & = & 1+\left\Vert F\right\Vert
_{1,q,p,\mathbf{\Theta }}+\left\Vert LF\right\Vert _{q-2,p,\mathbf{\Theta }}
\\
\mathrm{Q}_{F,\overline{F},\mathbf{\Theta }}(q,p) & = & 1+\left\Vert
F\right\Vert _{1,q,p,\mathbf{\Theta }}+\left\Vert LF\right\Vert _{q-2,p,%
\mathbf{\Theta }}+\left\Vert \overline{F}\right\Vert _{1,q,p,\mathbf{\Theta }%
}+\left\Vert L\overline{F}\right\Vert _{q-2,p,\mathbf{\Theta }}%
\end{array}
\label{unicadef}
\end{equation}%
with the convention ${\mathrm{S}}_{F,\mathbf{\Theta }}(p)=+\infty $ if the
r.h.s. is not finite.


\begin{proposition}
\label{prop-IBPU} Let $\kappa \in {\mathbb{N}}^{\ast }$ and assume that $%
\mathrm{m}_{\kappa ,p}(\mathbf{\Theta })<\infty $ for all $p\geq 1$.
Let $F\in \mathcal{S}^{d}$ be such that ${\mathrm{S}}_{F,\mathbf{\Theta }%
}(p)<\infty $ for every $p\in {\mathbb{N}}$. 
Let $\gamma _{F}$ be the inverse of $\sigma _{F}$ on the set $\{\mathbf{%
\Theta }\neq 0\}$. Then the following localized integration by parts formula
holds: for every $f\in C_{b}^{\infty }({\mathbb{R}}^{d})$, $G\in \mathcal{S}$
and for every multi index $\alpha $ of length equal to $q\leq \kappa $ one
has
\begin{equation*}
{\mathbb{E}}_{\mathbf{\Theta }}(\partial _{\alpha }f(F)\,G)={\mathbb{E}}_{%
\mathbf{\Theta }}(f(F)H_{\alpha ,\mathbf{\Theta }}^{q}(F,G))
\end{equation*}%
where as $r=1,\ldots ,d$%
\begin{eqnarray}
H_{r,\mathbf{\Theta }}(F,G) &=&\sum_{r^{\prime }=1}^{d}G\gamma
_{F}^{r^{\prime },r}LF^{r^{\prime }}-\langle D(G\gamma _{F}^{r^{\prime
},r}),DF^{r^{\prime }}\rangle _{J}-G\gamma _{F}^{r^{\prime },r}\langle D(\ln
\mathbf{\Theta }),DF^{r^{\prime }}\rangle _{J}  \label{Mall5} \\
&=&H_{r}(F,G)-G\sum_{r^{\prime }=1}^{d}\gamma _{F}^{r^{\prime },r}\langle
D(\ln \mathbf{\Theta }),DF^{r^{\prime }}\rangle _{J}  \notag
\end{eqnarray}%
and for a general multi index $\beta $ with $|\beta |=q$
\begin{equation*}
H_{\beta ,\mathbf{\Theta }}^{q}(F,G)=H_{\beta _{q},\mathbf{\Theta }}\big(%
F,H_{(\beta _{1},\ldots ,\beta _{q-1}),\mathbf{\Theta }}^{q-1}(F,G)\big).
\end{equation*}
\end{proposition}

\textbf{Proof.} For $|\beta |=1$, the integration by parts formula
immediately follows from the equality ${\mathbb{E}}_{\mathbf{\Theta }%
}(\partial _{i}f(F)G)={\mathbb{E}}(\partial _{i}f(F)G\mathbf{\Theta })={%
\mathbb{E}}(f(F)H_{i}(F,G\mathbf{\Theta }))$, so that $H_{i,\mathbf{\Theta }%
}(F,G)=\frac{1}{\mathbf{\Theta }}H_{i}(F,G\mathbf{\Theta })$, and this gives
the formula for $H_{i,\mathbf{\Theta }}(F,G)$. For higher order integration
by parts it suffices to iterate this procedure. $\square $\medskip

We give now estimates for the weights in the integration by parts formula.

\begin{proposition}
\label{prop-IBPU-est} Let $\kappa \in {\mathbb{N}}^{\ast }$ and $l\in {%
\mathbb{N}}$ be such that $\mathrm{m}_{l+\kappa +1,p}(\mathbf{\Theta }%
)<\infty $ for all $p\geq 1$. Let $F,\overline{F}\in \mathcal{S}^{d}$, with $%
{\mathrm{S}}_{F,\mathbf{\Theta }}(p),{\mathrm{S}}_{\overline{F},\mathbf{%
\Theta }}(p)<\infty $ for every $p$, and $G,\overline{G}\in \mathcal{S}$.
For $q\leq \kappa$, let $H_{\beta ,\mathbf{\Theta }}^{q}(\cdot ,\cdot )$ be
the weight of the integration by parts formula as in Proposition \ref%
{prop-IBPU}. Then for every $p\geq 1$ one may find two universal constants $%
C,p^{\prime }\in \mathcal{C(}\kappa ,d\mathcal{)}$ such that for every multi
index $\beta $ with $|\beta |=q\leq \kappa $
\begin{equation}
\Vert H_{\beta ,\mathbf{\Theta }}^{q}(F,G)\Vert _{l,p,\mathbf{\Theta }}\leq
C\,B_{l+q,p^{\prime },\mathbf{\Theta }}(F)^{q}\,\Vert G\Vert _{l+q,p^{\prime
},\mathbf{\Theta }},  \label{Mal21U}
\end{equation}%
and

\begin{equation}  \label{Mall22U}
\begin{array}{l}
\Vert H_{\beta ,\mathbf{\Theta }}^{q}(F,G)-H_{\beta ,\mathbf{\Theta }}^{q}(%
\overline{F},\overline{G})\Vert _{l,p,\mathbf{\Theta }} \leq \smallskip \\
\qquad\quad \leq C\,B_{l+q,p^{\prime },\mathbf{\Theta }}(F,\overline{F}%
)^{q(q+1)/2}\,
\big(1+\|G\|_{l+q,p'}+\|\overline{G}\|_{l+q,p'}\big)
\times \smallskip \\
\qquad\quad\quad \times (\Vert F-\overline{F}\Vert _{l+q+1,p,\mathbf{\Theta }%
}+\Vert LF-L\overline{F}\Vert _{l+q-1,p,\mathbf{\Theta }}+\Vert G-\overline{G%
}\Vert _{l+q,p,\mathbf{\Theta }})%
\end{array}%
\end{equation}
where
\begin{align}
& B_{l,p,\mathbf{\Theta }}(F)={\mathrm{S}}_{F,\mathbf{\Theta }}(p)^{l+1}%
\mathrm{Q}_{F,\mathbf{\Theta }}(l+1,p)^{2d(l+2)}\mathrm{m}_{l,p}(\mathbf{%
\Theta })  \label{Mall22bU} \\
& B_{l,p,\mathbf{\Theta }}(F,\overline{F})={\mathrm{S}}_{F,\mathbf{\Theta }%
}(p)^{l+1}{\mathrm{S}}_{\overline{F},\mathbf{\Theta }}(p)^{l+1}\mathrm{Q}_{F,%
\overline{F},\mathbf{\Theta }}(l+1,p)^{2d(l+2)}\mathrm{m}_{l,p}(\mathbf{%
\Theta })
\end{align}%
${\mathrm{S}}_{\cdot ,\mathbf{\Theta }}(p)$, $\mathrm{Q}_{\cdot ,\mathbf{%
\Theta }}(l,p)$ and $\mathrm{Q}_{\cdot ,\cdot ,\mathbf{\Theta }}(l,p)$ being
defined in \eqref{unicadef}.
\end{proposition}

\textbf{Proof.} By using the same arguments as in Theorem \ref{thHq}, one
gets that there exists $C\in\mathcal{C}(q,d)$ such that for every multi
index $\beta $ of length $q$ then
\begin{equation*}
|H_{\beta ,\mathbf{\Theta }}^{q}(F,G)|_{l}\leq C\,A_{l+q}(F)^{q}\big(1+|D\ln
\mathbf{\Theta}|_{l+q-1}\big)^{q}|G|_{l+q}
\end{equation*}%
and
\begin{align*}
|H_{\beta ,\mathbf{\Theta }}^{q}(F,G)-H_{\beta ,\mathbf{\Theta }}^{q}(%
\overline{F},\overline{G})|_{l}\leq & C\,A_{l+q}(F,\overline{F})^{\frac{%
q(q+1)}{2}}\big(1+|D\ln \mathbf{\Theta}|_{l+q-1}\big)^{\frac{q(q+1)}{2}}
\times \\
& \times \big(1+|G|_{l+q}+|\overline{G}|_{l+q}\big)^{q}\times \\
& \times \big(|F-\overline{F}|_{l+q+1}+|L(F-\overline{F})|_{l+q-1}+|G-%
\overline{G}|_{l+q}\big),
\end{align*}%
where $A_{l}(F)$ and $A_{l}(F,\overline{F})$ are defined in \eqref{gamma11'}
and \eqref{AAn} respectively (as usual, $|\cdot |_{0}\equiv |\cdot |$). By
using H\"{o}lder inequality one gets \eqref{Mal21U} and (\ref{Mall22U}). $%
\square $\medskip

In next Lemma we study properties of $H_{\beta ,\mathbf{\Theta }}^{q}(F,G)$
in the case $G$ is a special function of $F$. We denote with $B_{r}(0)$ the
ball with radius $r$ centered at $0$.

\begin{lemma}
\label{rem-ind} Let $\phi \in C_{b}^{\infty }({\mathbb{R}}^{d})$ be such
that $\mbox{\large \bf 1}_{B_{1}(0)}\leq \phi \leq \mbox{\large \bf 1}%
_{B_{2}(0)}$ and set $\phi _{x}(y)=\phi (x-y)$. For $l\in {\mathbb{N}}$, $%
p\geq 1$ and $F,\overline{F}\in \mathcal{S}^{d}$, one has
\begin{align}
& \Vert \phi _{x}(F)\Vert _{l,p,\mathbf{\Theta }}\leq C(1+\Vert F\Vert
_{1,l,2p,\mathbf{\Theta }})^{l}\P_{\mathbf{\Theta}}(F\in B_{2}(x)),  \label{est-phi1} \\
& \Vert \phi _{x}(F)-\phi _{x}(\overline{F})\Vert _{l,p,\mathbf{\Theta }%
}\leq C\Vert F-\overline{F}\Vert _{l,2p,\mathbf{\Theta }}(1+\Vert F\Vert
_{1,l,2p,\mathbf{\Theta }}+\Vert \overline{F}\Vert _{1,l,2p,\mathbf{\Theta }%
})^{l},  \label{est-phi2}
\end{align}%
in which $C\in \mathcal{C(}l,p,d\mathcal{)}$. Moreover for every $F\in
\mathcal{S}^{d}$ and $V\in \mathcal{S}$ one may find universal constants $%
C,a,p^{\prime }\in \mathcal{C(}q,l,p,d\mathcal{)}$ such that for every multi
index $\beta $ with $|\beta |=q$
\begin{equation}
\Vert H_{\beta ,\mathbf{\Theta }}^{q}(F,V\phi _{x}(F))\Vert _{l,p,\mathbf{%
\Theta }}\leq {\mathrm{S}}_{F,\mathbf{\Theta }}(p^{\prime })^{a}\mathrm{Q}%
_{F,\mathbf{\Theta }}(l+q+1,p^{\prime })^{a}\mathrm{m}_{l+q,p}(\mathbf{%
\Theta })^{a}{\mathbb{P}}_{\mathbf{\Theta }}(|F-x|<2)^{\frac{1}{p^{\prime }}%
}\times \Vert V\Vert _{l,p^{\prime },\mathbf{\Theta }},  \label{est-phi0}
\end{equation}%
where ${\mathrm{S}}_{F,\mathbf{\Theta }}(p)$, $\mathrm{Q}_{F,\mathbf{\Theta }%
}(l,p)$ and $\mathrm{m}_{l,p}(\mathbf{\Theta })$ are defined in %
\eqref{unicadef}.
\end{lemma}

\textbf{Proof.} We prove \eqref{est-phi2}, \eqref{est-phi1} following with
similar arguments. First, for a multi-index $\alpha $ with $|\alpha |=k$,
one has
\begin{equation*}
D_{\alpha }\phi _{x}(F)=\sum_{l=1}^{k}\phi _{x}^{(l)}(F)\!\!\!\sum_{\beta
_{1},\ldots ,\beta _{l}\in \mathcal{B}_{\alpha }}\!\!D_{\beta _{1}}F\cdots
D_{\beta _{l}}F
\end{equation*}%
where $C>0$ depends on $l,d$ only and \textquotedblleft $\beta _{1},\ldots
,\beta _{l}\in \mathcal{B}_{\alpha }$\textquotedblright\ means that $\beta
_{1},\ldots ,\beta _{l}$ are non empty multi indexes of $\alpha $ running
through the list of all of the (non empty) \textquotedblleft
blocks\textquotedblright\ of $\alpha $. So, straightforward computations
give
\begin{align*}
& D_{\alpha }\phi _{x}(F)-D_{\alpha }\phi _{x}(\overline{F}) \\
& \quad =\sum_{l=1}^{k}\big(\phi _{x}^{(l)}(F)-\phi _{x}^{(l)}(\overline{F})%
\big)\!\!\!\sum_{\beta _{1},\ldots ,\beta _{l}\in \mathcal{B}_{\alpha
}}\!\!D_{\beta _{1}}F\cdots D_{\beta _{l}}F+ \\
& \qquad +\sum_{l=1}^{k}\phi _{x}^{(l)}(\overline{F})\!\!\!\sum_{\beta
_{1},\ldots ,\beta _{l}\in \mathcal{B}_{\alpha }}\!\!\sum_{j=1}^{l}\Big(%
\prod_{k=1}^{j-1}D_{\beta _{k}}F\Big)\big(D_{\beta _{j}}F-D_{\beta _{j}}%
\overline{F}\big)\Big(\prod_{k=j+1}^{l}D_{\beta _{k}}\overline{F}\Big)
\end{align*}%
with the understanding $\prod_{k=1}^{0}(\cdot
)_{k}=1=\prod_{k=l+1}^{l}(\cdot )_{k}$. Since $\phi _{x}^{(l)}$ is Lipschitz
continuous, with a Lipschitz constant independent of $x$, it follows that
\begin{align*}
|\phi _{x}(F)-\phi _{x}(\overline{F})|_{l}& \leq C|F-\overline{F}%
|(1+|F|_{1,l})^{l}+C|F-\overline{F}|_{l}(1+|F|_{1,l}+|\overline{F}%
|_{1,l})^{l-1} \\
& \leq C|F-\overline{F}|_{l}(1+|F|_{1,l}+|\overline{F}|_{1,l})^{l}
\end{align*}%
and by using the H\"{o}lder inequality one gets \eqref{est-phi2}.

\smallskip

As for \eqref{est-phi0}, we first note that since $\phi _{x}(y)\equiv 0$ for
$|y-x|>2$ then \eqref{Mall5} gives
\begin{equation}
D_{\alpha }H_{\beta ,\mathbf{\Theta }}^{q}(F,V\phi _{x}(F))=D_{\alpha
}H_{\beta ,\mathbf{\Theta }}^{q}(F,V\phi _{x}(F))\mbox{\large \bf 1}%
_{\{|F-x|<2\}}  \label{est-phi}
\end{equation}%
for every multi index $\alpha $. So, for $l\in {\mathbb{N}}$ we can write
\begin{equation*}
|H_{\beta ,\mathbf{\Theta }}^{q}(F,V\phi _{x}(F))|_{l}=|H_{\beta ,\mathbf{%
\Theta }}^{q}(F,V\phi _{x}(F))|_{l}\mbox{\large \bf
1}_{\{|F-x|<2\}}.
\end{equation*}%
Therefore, \eqref{est-phi0} is a consequence of the use of the H\"{o}lder
inequality and of the estimate \eqref{Mal21U}. $\square $\medskip

We recall that the Poisson kernel $Q_{d}$ is the solution to the equation $%
\Delta Q_{d}=\delta _{0}$ in ${\mathbb{R}}^{d}$ ($\delta_0$ denoting the
Dirac mass in $\{0\}$) and has the following explicit form:
\begin{equation}
Q_1(x)=\max\{x,0\},\quad Q_{2}(x)=a_{2}^{-1}\ln \left\vert x\right\vert
\quad \mbox{and}\quad Q_{d}(x)=-a_{d}^{-1}\left\vert x\right\vert ^{2-d},d>2
\label{den4}
\end{equation}%
where $a_{d}$ is the area of the unit sphere in ${\mathbb{R}}^{d}.$ By using
the result in \cite{bib:[BC0]}, we have the following

\begin{proposition}
\label{1} Let $\phi \in C_{b}^{\infty }({\mathbb{R}}^{d})$ be such that $%
1_{B_{1}(0)}\leq \phi \leq 1_{B_{2}(0)}$ and set $\phi _{x}(y)=\phi (x-y)$.
Let $\kappa \in {\mathbb{N}}^{\ast }$ and assume that $\mathrm{m}_{\kappa,p}(%
\mathbf{\Theta })<\infty $ for every $p\geq 1$. Let $F\in \mathcal{S}^{d}$
be such that ${\mathrm{S}}_{F,\mathbf{\Theta }}(p)<\infty $ for every $p\geq
1$.

\medskip

\textbf{A}. Let $Q_{d}$ be the Poisson kernel in ${\mathbb{R}}^{d}$ given in
(\ref{den4}). Then for every $p>d$ there exists a universal constant $C\in%
\mathcal{C}(d,p) $ such that
%
%
\begin{equation}
\Vert \nabla Q_{d}(F-x)\Vert _{\frac{p}{p-1},U}\leq C\Vert H_{\mathbf{\Theta
}}(F,1)\Vert _{p,\mathbf{\Theta }}^{k_{p,d}}  \label{Mall5'}
\end{equation}%
%
%
%
%
%
%
%
%
%
%
%
%
%
%
%
where $k_{p,d}=(d-1)/(1-d/p)$ and $H_{\mathbf{\Theta }}(F,1)$ denotes the
vector in ${\mathbb{R}}^{d}$ whose $i$th entry is given be $H_{i,\mathbf{%
\Theta }}(F,1)$.

\smallskip

\textbf{B}. Under ${\mathbb{P}}_{\mathbf{\Theta }}$, the law of $F$ is
absolutely continuous and has a density $p_{F,\mathbf{\Theta }}\in C^{\kappa
-1}({\mathbb{R}}^{d})$ whose derivatives up to order $\kappa -1$ may be
represented as%

\begin{equation}
\partial _{\alpha }p_{F,\mathbf{\Theta }}(x)=\sum_{i=1}^{d}{\mathbb{E}}_{%
\mathbf{\Theta }}\big(\partial _{i}Q_{d}(F-x)H_{(i,\alpha ),\mathbf{\Theta }%
}^{q+1}(F,\phi _{x}(F))\big)  \label{Mall5''}
\end{equation}%
for every multi index $\alpha $ with $|\alpha |=q\leq \kappa-1 $.

\smallskip

\textbf{C.} Let $V$ a random variable taking values in $(0,1)$ and such that $\mathrm{m}_{\kappa
,p}(V)<\infty $ for every $p\geq 1$. Then for $|\alpha |=q\leq \kappa -1$
one has
\begin{equation}
|\partial _{\alpha }p_{F,\mathbf{\Theta }V}(x)|\leq C{\mathrm{S}}_{F,\mathbf{%
\Theta }}(p^{\prime })^{a}\mathrm{Q}_{F,\mathbf{\Theta }}(q+2,p^{^{\prime
}})^{a}\mathrm{m}_{q+1,p^{\prime }}(\Theta )^{a}\left\Vert V\right\Vert
_{q+1,p^{\prime },\mathbf{\Theta }}\times {\mathbb{P}}_{\mathbf{\Theta}}(|F-x|<2)^{b},
\label{Bis3}
\end{equation}%
in which $C,a,b,p^{\prime }\in \mathcal{C(}\kappa ,d\mathcal{)}$.
\end{proposition}

\textbf{Proof}. \textbf{A.} This point is actually Theorem 5 in \cite%
{bib:[BC0]} 
(recall that $\Vert 1\Vert _{W_{\mu _{F}}^{1,p}}\leq \Vert H(F,1)\Vert _{p}$%
, see Remark 17 in \cite{bib:[BC0]}) with ${\mathbb{P}}$ replaced by ${%
\mathbb{P}}_{\mathbf{\Theta }}$.

\smallskip

\textbf{B.} Set $\mu _{F,\mathbf{\Theta }}$ the law of $F$ under ${\mathbb{P}%
}_{\mathbf{\Theta }}$ and let $\alpha $ denote a multi index with $|\alpha
|=q$. By using the arguments similar to the ones developed in Proposition 10
in \cite{bib:[BC0]} one easily gets (notations from that paper)
\begin{equation*}
\partial _{\alpha }p_{F,\mathbf{\Theta }}(x)=(-1)^{|\alpha
|+1}\sum_{i=1}^{d}\int_{{\mathbb{R}}^{d}}\partial _{i}Q_{d}(y-x)\partial
_{(i,\alpha )}^{\mu _{F,\mathbf{\Theta }}}\phi _{x}(y)\mu _{F,\mathbf{\Theta
}}(dy).
\end{equation*}%
And by recalling that $(-1)^{|\alpha |+1}\partial _{(i,\alpha )}^{\mu _{F,%
\mathbf{\Theta}}}\phi _{x}(F)={\mathbb{E}}_{\mathbf{\Theta }}(H_{(i,\alpha ),%
\mathbf{\Theta }}^{q+1}(F,\phi _{x}(F))\mid F)$ (see Section 3 of \cite%
{bib:[BC0]}), \eqref{Mall5''} follows.

\textbf{C.} We first note that $\mathrm{m}_{\mathbf{\Theta }V}(\kappa
,p)\leq C(\mathrm{m}_{\mathbf{\Theta }}(\kappa ,p)+\mathrm{m}_{V}(\kappa
,p)) $. So, we can apply \eqref{Mall5''} with localization $\mathbf{\Theta }%
V $ and we get
\begin{equation*}
\partial _{\alpha }p_{F,\mathbf{\Theta }V}(x)=\sum_{i=1}^{d}{\mathbb{E}}_{%
\mathbf{\Theta}}\big(\partial _{i}Q_{d}(F-x)VH_{(i,\alpha ),\mathbf{\Theta }%
V}^{q+1}(F,\phi _{x}(F))\big).
\end{equation*}%
Now, from \eqref{Mall5} one has $VH_{i,\mathbf{\Theta }V}(F,\phi
_{x}(F))=H_{i,\mathbf{\Theta }}(F,V\phi _{x}(F))$ and by iteration it
follows that $VH_{\beta ,\mathbf{\Theta }V}^{q}(F,\phi _{x}(F))=H_{\beta ,%
\mathbf{\Theta }}^{q}(F,V\phi _{x}(F))$. Therefore,
\begin{equation*}
\partial _{\alpha }p_{F,\mathbf{\Theta }V}(x)=\sum_{i=1}^{d}{\mathbb{E}}_{%
\mathbf{\Theta}}\big(\partial _{i}Q_{d}(F-x)H_{(i,\alpha ),\mathbf{\Theta }%
}^{q+1}(F,V\phi _{x}(F))\big)
\end{equation*}%
and, by using the H\"{o}lder inequality, for $p>d$ we have
\begin{align*}
|\partial _{\alpha }p_{F,\mathbf{\Theta }V}(x)|& \leq \sum_{i=1}^{d}\Vert
\nabla Q_{d}(F-x)\Vert _{\frac{p}{p-1},\mathbf{\Theta }}\Vert H_{(i,\alpha ),%
\mathbf{\Theta }}^{q+1}(F,V\phi _{x}(F))\Vert _{p,\mathbf{\Theta }} \\
& \leq \sum_{i=1}^{d}\Vert H_{\mathbf{\Theta }}(F,1)\Vert _{p,\mathbf{\Theta
}}^{k_{p,d}}\Vert H_{(i,\alpha ),\mathbf{\Theta }}^{q+1}(F,V\phi
_{x}(F))\Vert _{p,\mathbf{\Theta }}
\end{align*}%
in which we have used \eqref{Mall5'}. Now, by using \eqref{Mal21U} to
estimate the first term and by applying \eqref{est-phi0} to the second one, %
\eqref{Bis3} follows. $\square $

\subsection{The distance between density functions and their derivatives}

\label{sect-dist}

We compare now the probability density functions (and their derivatives) of
two random variables under ${\mathbb{P}}_{\mathbf{\Theta }}.$

\begin{proposition}
\label{prop-dist} Let $q\in {\mathbb{N}}$ and assume that $\mathrm{m}%
_{q+2,p}(\mathbf{\Theta })<\infty $ for every $p\geq 1$. Let $F,G\in
\mathcal{S}^{d}$ be such that
\begin{equation}
{\mathrm{S}}_{F,G,\mathbf{\Theta }}(p):=1+\sup_{0\leq \lambda \leq 1}\Vert
(\det \sigma _{G+\lambda (F-G)})^{-1}\Vert _{p,\mathbf{\Theta }}<\infty
,\quad \forall p\in {\mathbb{N}}.  \label{mFGU}
\end{equation}%
Then under ${\mathbb{P}}_{\mathbf{\Theta }}$ the laws of $F$ and $G$ are
absolutely continuous with respect to the Lebesgue measure with density $%
p_{F,\mathbf{\Theta }}$ and $p_{G,\mathbf{\Theta }}$ respectively and for
every multi index $\alpha $ with $|\alpha |=q$ there exist constants $%
C,a,b,p^{\prime }\in \mathcal{C(}q,d\mathcal{)}$ such that%
\begin{equation}  \label{Mall2}
\begin{array}{rcl}
\left\vert \partial _{\alpha }p_{F,\mathbf{\Theta }}(y)-\partial _{\alpha
}p_{G,\mathbf{\Theta }}(y)\right\vert & \leq & C{\mathrm{S}}_{F,G,\mathbf{%
\Theta }}(p^{\prime })^{a}\mathrm{Q}_{F,G,\mathbf{\Theta }}(q+3,p^{\prime
})^{a}\mathrm{m}_{q+2,p}^{a}(\mathbf{\Theta })\times \smallskip \\
&  & \times (\Vert F-G\Vert _{q+2,p^{\prime },\mathbf{\Theta }}+\Vert
LF-LG\Vert _{q,p^{\prime },\mathbf{\Theta }}) \times \smallskip \\
&  & \times ({\mathbb{P}}_\mathbf{\Theta}(\left\vert F-y\right\vert < 2)+{\mathbb{P}}_\mathbf{\Theta}
(\left\vert G-y\right\vert < 2))^{b}%
\end{array}%
\end{equation}
with $\mathrm{m}_{k,p}(\mathbf{\Theta })$ and $\mathrm{Q}_{F,G,\mathbf{%
\Theta }}(k,p)$ given in (\ref{Mall1}) and \eqref{unicadef} respectively.
\end{proposition}

\textbf{Proof}. Throughout this proof, $C,p^{\prime },a,b\in \mathcal{C(}q,d%
\mathcal{)}$ will denote constants that can vary from line to line. By
applying Lemma \ref{1}, under ${\mathbb{P}}_{\mathbf{\Theta }}$ the laws of $%
F$ and $G$ are both absolutely continuous with respect to the Lebesgue
measure and for every multi index $\alpha $ with $|\alpha |=q$ one has
\begin{align*}
\partial _{\alpha }p_{F,\mathbf{\Theta }}(y)-\partial _{\alpha }p_{G,\mathbf{%
\Theta }}(y)=& \sum_{j=1}^{d}{\mathbb{E}}_{\mathbf{\Theta }}\big(\big(%
\partial _{j}Q_{d}(F-y)-\partial _{j}Q_{d}(G-y)\big)H_{(j,\alpha ),\mathbf{%
\Theta }}^{q+1}(G,\phi _{y}(G))\big)+ \\
& +\sum_{j=1}^{d}{\mathbb{E}}_{\mathbf{\Theta }}\big(\partial _{j}Q_{d}(F-y)%
\big(H_{(j,\alpha ),\mathbf{\Theta }}^{q+1}(F,\phi _{y}(F))-H_{(j,\alpha ),%
\mathbf{\Theta }}^{q+1}(G,\phi _{y}(G))\big)\big) \\
=& :\sum_{j=1}^{d}I_{j}+\sum_{j=1}^{d}J_{j}.
\end{align*}%
By using (\ref{Mall5'}), for $p>d$ we obtain%
\begin{align*}
|J_{j}|& \leq C\,\Vert \nabla Q_{d}(F-y)\Vert _{\frac{p}{p-1},U}\Vert
H_{(j,\alpha ),\mathbf{\Theta }}^{q+1}(F,\phi _{y}(F))-H_{(j,\alpha ),%
\mathbf{\Theta }}^{q+1}(G,\phi _{y}(G))\Vert _{p,\mathbf{\Theta }} \\
& \leq C\,\Vert H_{\mathbf{\Theta }}(F,1)\Vert _{p,\mathbf{\Theta }%
}^{k_{d,p}}\Vert H_{(j,\alpha ),\mathbf{\Theta }}^{q+1}(F,\phi
_{y}(F))-H_{(j,\alpha ),\mathbf{\Theta }}^{q+1}(G,\phi _{y}(G))\Vert _{p,%
\mathbf{\Theta }}.
\end{align*}%
Now, from \eqref{est-phi} (with $\alpha =\emptyset $) it follows that the
above term is null on $\{\left\vert F-y\right\vert \geq 2\}\cap \{\left\vert
G-y\right\vert \geq 2\}.$ So%
\begin{align*}
& |H_{(j,\alpha ),\mathbf{\Theta }}^{q+1}(F,\phi _{y}(F))-H_{(j,\alpha ),%
\mathbf{\Theta }}^{q+1}(G,\phi _{y}(G))| \\
& \quad \leq |H_{(j,\alpha ),\mathbf{\Theta }}^{q+1}(F,\phi
_{y}(F))-H_{(j,\alpha ),\mathbf{\Theta }}^{q+1}(G,\phi _{y}(G))|%
\mbox{\large
\bf 1}_{\{|F-y|<2\}}+ \\
& \qquad +|H_{(j,\alpha ),\mathbf{\Theta }}^{q+1}(F,\phi
_{y}(F))-H_{(j,\alpha ),\mathbf{\Theta }}^{q+1}(G,\phi _{y}(G))|%
\mbox{\large
\bf 1}_{\{|G-y|<2\}}
\end{align*}%
so that the H\"{o}lder inequality gives
\begin{align*}
|J_{j}|\leq & C\,\Vert H_{\mathbf{\Theta }}(F,1)\Vert _{p,U}^{k_{d,p}}\Vert
H_{(j,\alpha ),\mathbf{\Theta }}^{q+1}(F,\phi _{y}(F))-H_{(j,\alpha ),%
\mathbf{\Theta }}^{q+1}(G,\phi _{y}(G))\Vert _{2p,U}\times  \\
& \times \big({\mathbb{P}}_{\mathbf{\Theta }}(|F-x|<2)+{\mathbb{P}}_{\mathbf{%
\Theta }}(|G-x|<2)\big)^{\frac{1}{2p}}.
\end{align*}%
So, by applying \eqref{Mal21U} and (\ref{Mall22U}), there exists $p^{\prime
}>p>d$ such that
\begin{align*}
|J_{j}|\leq & CB_{q+1,p^{\prime },\mathbf{\Theta }}(F,G)^{k_{d,p}+\frac{%
(q+1)(q+2)}{2}}\times  \\
& \times \big(\Vert F-G\Vert _{q+2,p^{\prime },\mathbf{\Theta }}+\Vert
L(F-G)\Vert _{q,p^{\prime },\mathbf{\Theta }}+\Vert \phi _{y}(F)-\phi
_{y}(G)\Vert _{q+1,p^{\prime },\mathbf{\Theta }}\big)\times  \\
& \times \big({\mathbb{P}}_{\mathbf{\Theta }}(|F-x|<2)+{\mathbb{P}}_{\mathbf{%
\Theta }}(|G-x|<2)\big)^{\frac{1}{2p^{\prime }}},
\end{align*}%
$B_{q+1,p^{\prime },\mathbf{\Theta }}(F,G)$ being defined in \eqref{Mall22bU}%
. By using \eqref{est-phi2} and the quantities ${\mathrm{S}}_{F,G,\mathbf{%
\Theta }}(p)$ and $\mathrm{Q}_{F,G,\mathbf{\Theta }}(k,p)$, for a suitable $%
a>1$ and $p^{\prime }>d$ we can write
\begin{align*}
|J_{j}|\leq & C{\mathrm{S}}_{F,G,\mathbf{\Theta }}(p^{\prime })^{a}\mathrm{Q}%
_{F,G,\mathbf{\Theta }}(q+2,p^{\prime })^{a}\mathrm{m}_{q+1,p^{\prime }}(%
\mathbf{\Theta })^{a}\times \big(\Vert F-G\Vert _{q+2,p^{\prime },\mathbf{%
\Theta }}+\Vert L(F-G)\Vert _{q,p^{\prime },\mathbf{\Theta }}\big)\times  \\
& \times \big({\mathbb{P}}_{\mathbf{\Theta }}(|F-x|<2)+{\mathbb{P}}_{\mathbf{%
\Theta }}(|G-x|<2)\big)^{\frac{1}{2p^{\prime }}}.
\end{align*}%
We study now $I_{j}$. For $\lambda \in \lbrack 0,1]$ we denote $F_{\lambda
}=G+\lambda (F-G)$ and we use Taylor's expansion to obtain%
\begin{equation*}
I_{j}=\sum_{k=1}^{d}R_{k,j}\quad \mbox{with}\quad R_{k,j}=\int_{0}^{1}{%
\mathbb{E}}_{\mathbf{\Theta }}\big(\partial _{k}\partial
_{j}Q_{d}(F_{\lambda }-y)H_{(j,\alpha ),\mathbf{\Theta }}^{q+1}(G,\phi
_{y}(G))(F-G)_{k}\big)d\lambda .
\end{equation*}%
Let $V_{k,j}=H_{(j,\alpha ),\mathbf{\Theta }}^{q+1}(G,\phi _{y}(G))(F-G)_{k}.
$ Since for $\lambda \in \lbrack 0,1]$ then ${\mathbb{E}}_{\mathbf{\Theta }%
}((\det \sigma _{F_{\lambda }})^{-p}))$ $<\infty $ for every $p$, we can use
the integration by parts formula with respect to $F_{\lambda }$, so%
\begin{equation*}
R_{k,j}=\int_{0}^{1}{\mathbb{E}}_{\mathbf{\Theta }}\big(\partial
_{j}Q_{d}(F_{\lambda }-y)H_{k,\mathbf{\Theta }}(F_{\lambda },V_{k,j})\big)%
d\lambda .
\end{equation*}%
Therefore, by taking $p>d$ and by using again (\ref{Mall5'}), \eqref{Mal21U}
and (\ref{Mall22U}), we get
\begin{align*}
|R_{k,j}|\leq & \int_{0}^{1}\Vert \partial _{j}Q_{d}(F_{\lambda }-y)\Vert _{%
\frac{p}{p-1},\mathbf{\Theta }}\Vert H_{k,\mathbf{\Theta }}(F_{\lambda
},V_{k,j})\Vert _{p,\mathbf{\Theta }}d\lambda  \\
\leq & C\int_{0}^{1}\Vert H_{\mathbf{\Theta }}(F_{\lambda },1)\Vert
_{p,U}^{k_{d,p}}\Vert H_{k,\mathbf{\Theta }}(F_{\lambda },V_{k,j})\Vert _{p,%
\mathbf{\Theta }}d\lambda  \\
\leq & C\int_{0}^{1}B_{1,p^{\prime },\mathbf{\Theta }}(F_{\lambda
})^{k_{d,p}+1}\Vert V_{k,j}\Vert _{1,p^{\prime }\mathbf{\Theta }}\,d\lambda
\end{align*}%
in which we have used \eqref{Mal21U}. Now, from \eqref{mFGU} and (\ref%
{Mall22bU}) it follows that
\begin{equation*}
B_{1,p,\mathbf{\Theta }}(F_{\lambda })\leq C{\mathrm{S}}_{F,G,\mathbf{\Theta
}}(p)^{4}\mathrm{Q}_{F,G,\mathbf{\Theta }}(2,p)^{8d}\times \mathrm{m}_{1,p}(%
\mathbf{\Theta }).
\end{equation*}%
Moreover,
\begin{equation*}
\Vert V_{k,j}\Vert _{1,p,\mathbf{\Theta }}=\Vert H_{(j,\alpha ),\mathbf{%
\Theta }}^{q+1}(G,\phi _{y}(G))(F-G)_{k}\Vert _{1,p,\mathbf{\Theta }}\leq
\Vert H_{(j,\alpha ),\mathbf{\Theta }}^{q+1}(G,\phi _{y}(G))\Vert _{1,2p,%
\mathbf{\Theta }}\Vert F-G\Vert _{1,2p,\mathbf{\Theta }}
\end{equation*}%
and from \eqref{est-phi0} and \eqref{Mal21U} we get
\begin{equation*}
\Vert V_{k,j}\Vert _{1,p,\mathbf{\Theta }}\leq C\,B_{q+2,p^{\prime },\mathbf{%
\Theta }}(G)^{q+1}\,\Vert \phi _{y}(G)\Vert _{q+2,p^{\prime },U}{\mathbb{P}}%
_{\mathbf{\Theta }}(|G-y|<2)^{\frac{1}{p^{\prime }}}\Vert F-G\Vert
_{1,p^{\prime },\mathbf{\Theta }}
\end{equation*}

By using (\ref{Mall22bU}),
\begin{equation*}
B_{q+2,p^{\prime },\mathbf{\Theta }}(G)\leq {\mathrm{S}}_{F,G,\mathbf{\Theta
}}(p)^{2q+4}\mathrm{Q}_{G,\mathbf{\Theta }}(q+3,p^{\prime})^{2d(q+3)}\times
\mathrm{m}_{q+2,p^{\prime }}(\mathbf{\Theta }),
\end{equation*}%
$\mathrm{Q}_{G,\mathbf{\Theta }}(l,p)$ being given in \eqref{unicadef}. We
use also \eqref{est-phi1} and, by inserting everything, we can resume by
writing
\begin{equation*}
|I_{j}|\leq C{\mathrm{S}}_{F,G,\mathbf{\Theta }}(p^{\prime})^{ a}\mathrm{Q}%
_{F,G,\mathbf{\Theta }}(q+3,p^{\prime })^{a}\mathrm{m}_{q+2,p^{\prime }}(%
\mathbf{\Theta })^{a}\Vert F-G\Vert _{1,p^{\prime },\mathbf{\Theta }}\times {%
\mathbb{P}}_{\mathbf{\Theta }}(|G-x|<2)^{b}
\end{equation*}%
and the statement follows. $\square $

\medskip

Using the localizing function in (\ref{Mall10}) and by applying Proposition %
\ref{prop-dist} we get the following result.

\begin{theorem}
\label{Diference} Let $q\in {\mathbb{N}}$. Assume that $\mathrm{m}_{q+2,p}(%
\mathbf{\Theta })<\infty $ for every $p\geq 1$. Let $F,G\in \mathcal{S}^{d}$
be such that ${\mathrm{S}}_{F,\mathbf{\Theta }}(p),{\mathrm{S}}_{G,\mathbf{%
\Theta }}(p)<\infty $ for every $p\in {\mathbb{N}}$. Then under ${\mathbb{P}}%
_{\mathbf{\Theta }}$, the laws of $F$ and $G$ are absolutely continuous with
respect to the Lebesgue measure, with densities $p_{F,\mathbf{\Theta }}$ and
$p_{G,\mathbf{\Theta }}$ respectively. Moreover, there exist constants $%
C,a,b,p^{\prime }\in \mathcal{C}(q,d)$ such that for every multi index $%
\alpha $ of length $q$ one has%
\begin{equation}
\begin{array}{rcl}
\left\vert \partial _{\alpha }p_{F,\mathbf{\Theta }}(y)-\partial _{\alpha
}p_{G,\mathbf{\Theta }}(y)\right\vert  & \leq  & C{\mathrm{S}}_{F,\mathbf{%
\Theta }}(p^{\prime })^{a}{\mathrm{S}}_{G,\mathbf{\Theta }}(p^{\prime })^{a}%
\mathrm{Q}_{F,G,\mathbf{\Theta }}(q+3,p^{\prime })^{a}\mathrm{m}_{q+2,p}^{a}(%
\mathbf{\Theta })\times \smallskip  \\
&  & \times (\Vert F-G\Vert _{q+2,p^{\prime },\mathbf{\Theta }}+\Vert
LF-LG\Vert _{q,p^{\prime },\mathbf{\Theta }})\times \smallskip  \\
&  & \times ({\mathbb{P}}_\mathbf{\Theta}(\left\vert F-y\right\vert <2)+{\mathbb{P}}_\mathbf{\Theta}
(\left\vert G-y\right\vert <2))^{b}%
\end{array}
\label{Mall12bis}
\end{equation}%
with $\mathrm{m}_{k,p}(\mathbf{\Theta })$ and $\mathrm{Q}_{F,G,\mathbf{%
\Theta }}(k,p)$ given in (\ref{Mall1}) and \eqref{unicadef} respectively.
\end{theorem}

\textbf{Proof}. The proof consists in proving that (\ref{Mall2}) holds.

Set $R=F-G.$ We use the deterministic estimate (\ref{b}) on the distance
between the determinants of two Malliavin covariance matrices: for every $%
\lambda \in \lbrack 0,1]$ we can write
\begin{align*}
\left\vert \det \sigma _{G+\lambda R}-\det \sigma _{G}\right\vert & \leq
C_{d}\left\vert DR\right\vert (\left\vert DG\right\vert +\left\vert
DF\right\vert )^{2d-1} \\
& \leq \big(\alpha _{d}|DR|^{2}(|DG|^{2}+|DF|^{2})^{\frac{2d-1}{2}}\big)%
^{1/2},
\end{align*}%
so that
\begin{equation}
\det \sigma _{G+\lambda R}\geq \det \sigma _{G}-\alpha _{d}\big(%
|DR|^{2}(|DG|^{2}+|DF|^{2})^{\frac{2d-1}{2}}\big)^{1/2}.  \label{alpha-d}
\end{equation}%
For $\psi _{a}$ as in (\ref{Mall10}), we define
\begin{equation*}
V=\psi _{1/8}(H)\quad \mbox{with}\quad H=|DR|^{2}\,\frac{%
(|DG|^{2}+|DF|^{2})^{\frac{2d-1}{2}}}{(\det \sigma _{G})^{2}},
\end{equation*}%
so that
\begin{equation}
V\neq 0\quad \Rightarrow \quad \det \sigma _{G+\lambda R}\geq \frac{1}{2}%
\det \sigma _{G}.  \label{V12}
\end{equation}%
Before continuing, let us give the following estimate for the Sobolev norm
of $H$. First, coming back to the notation $|\cdot |_{l}$ as in (\ref{norm1}%
), by using \eqref{prod} one easily get
\begin{equation*}
|H|_{l}\leq C|(\det \sigma _{G})^{-1}|_{l}^{2}\big(1+\big||DF|^{2}\big|_{l}+%
\big||DG|^{2}\big|_{l}\big)^{d}\big||DR|^{2}\big|_{l}.
\end{equation*}%
By using the estimate concerning the determinant from \eqref{est-det} and
the straightforward estimate $\big||DF|^{2}\big|_{l}\leq C|F|_{1,l+1}^{l+1}$,
we have
\begin{equation*}
|H|_{l}\leq C|(\det \sigma _{G})^{-1}|^{2(l+1)}\big(1+|F|_{1,l+1}+|G|_{1,l+1}%
\big)^{8dl}\,|R|_{ 1,  l+1}^{l+1}.
\end{equation*}%
As a consequence, by using the H\"{o}lder inequality we obtain
\begin{equation}
\Vert H\Vert _{l,p,\mathbf{\Theta }}\leq C{\mathrm{S}}_{G,\mathbf{\Theta }}(%
\bar{p})^{\bar{a}}\mathrm{Q}_{F,G,\mathbf{\Theta }}(l+1,\bar{p})^{\bar{a}%
}\Vert F-G\Vert _{l+1,\bar{p},\mathbf{\Theta }}^{l+1}  \label{est-Hl}
\end{equation}%
where $C,\bar{p},\bar{a}$ depends on $l,d,p$.

Now, because of \eqref{V12}, we have ${\mathrm{S}}_{F,G,\mathbf{\Theta }%
V}(p)\leq C{\mathrm{S}}_{G,\mathbf{\Theta }}(p)$, $C$ denoting a suitable
positive constant (which will vary in the following lines). We also have $%
\mathrm{m}_{k,p}(\mathbf{\Theta }V)\leq C(\mathrm{m}_{k,p}(\mathbf{\Theta })+%
\mathrm{m}_{k,p}(V))$. By (\ref{Mall11}) and \eqref{est-Hl} we have%
\begin{equation*}
\mathrm{m}_{q+2,p}(V)\leq C\,{\mathrm{S}}_{G,\mathbf{\Theta }}(\bar{p})^{%
\bar{a}}\mathrm{Q}_{F,G,\mathbf{\Theta }}(q+3,\bar{p})^{\bar{a}}
\end{equation*}%
for some $\bar{p},\bar{a}$, so that $m_{q+2,p}(\mathbf{\Theta }V)\leq C\,{%
\mathrm{S}}_{G,\mathbf{\Theta }}(\bar{p})^{\bar{a}}\mathrm{Q}_{F,G,\mathbf{%
\Theta }}(q+3,\bar{p})^{\bar{a}}\mathrm{m}_{q+2,\bar{p}}(\mathbf{\Theta })^{%
\bar{a}}$. So, we can apply (\ref{Mall2}) with localization $\mathbf{\Theta }%
V$ and we get%
\begin{equation*}
\begin{array}{rcl}
\left\vert \partial _{\alpha }p_{F,\mathbf{\Theta V}}(y)-\partial _{\alpha
}p_{G,\mathbf{\Theta V}}(y)\right\vert & \leq & C{\mathrm{S}}_{G,\mathbf{%
\Theta }}(p^{\prime })^{a}\mathrm{Q}_{F,G,\mathbf{\Theta }}(q+3,p^{\prime
})^{a}\mathrm{m}_{q+2,p}^{a}(\mathbf{\Theta })\times \smallskip \\
&  & \times (\Vert F-G\Vert _{q+2,p^{\prime },\mathbf{\Theta }}+\Vert
LF-LG\Vert _{q,p^{\prime },\mathbf{\Theta }})\times \smallskip \\
&  & \times ({\mathbb{P}}_\mathbf{\Theta}(\left\vert F-y\right\vert \leq 2)+{\mathbb{P}}_\mathbf{\Theta}
(\left\vert G-y\right\vert \leq 2))^{b} \notag%
\end{array}%
\end{equation*}
with $p^{\prime }>d$ and $C,a,b>0$ depending on $q,d$. We write now%
\begin{align*}
|\partial _{\alpha }p_{F,\mathbf{\Theta }}(y)-\partial _{\alpha }p_{G,%
\mathbf{\Theta }}(y)|\leq & \left\vert \partial _{\alpha }p_{F,\mathbf{%
\Theta }V}(y)-\partial _{\alpha }p_{G,\mathbf{\Theta }V}(y)\right\vert + \\
& +\left\vert \partial _{\alpha }p_{F,\mathbf{\Theta }(1-V)}(y)\right\vert
+\left\vert \partial _{\alpha }p_{G,\mathbf{\Theta }V}(y)\right\vert ,
\end{align*}%
and we have already seen that the first addendum on the r.h.s. behaves as
desired. So, it suffices to show that also the remaining two terms have the
right behavior. To this purpose, we use (\ref{Bis3}). We have
\begin{equation*}
|\partial _{\alpha }p_{F,\mathbf{\Theta }(1-V)}(x)|\leq C{\mathrm{S}}_{F,%
\mathbf{\Theta }}(p)^{a}\mathrm{Q}_{F,\mathbf{\Theta }}(q+2,p)^{a}\mathrm{m}_{q+1,p}(%
\mathbf{\Theta })^{a}\left\Vert 1-V\right\Vert _{q+1,p,\mathbf{\Theta }%
}\times {\mathbb{P}}_{\mathbf{\Theta }}(|F-x|<2)^{b}.
\end{equation*}%
Now, we can write
\begin{equation*}
\Vert 1-V\Vert _{q+1,p,\mathbf{\Theta }}^{p}={\mathbb{E}}_{\mathbf{\Theta }%
}(|1-V|^{p})+\Vert DV\Vert _{q,p,\mathbf{\Theta }}.
\end{equation*}%
But $1-V\neq 0$ implies that $H\geq 1/8$. Moreover, from (\ref{loc9}), $%
\Vert DV\Vert _{q,p,\mathbf{\Theta }}\leq \Vert V\Vert _{q+1,p,\mathbf{%
\Theta }}$ $\leq C\Vert H\Vert _{q+1,p(q+1),\mathbf{\Theta }}^{q+1}$. So, we
have
\begin{align*}
\Vert 1-V\Vert _{q+1,p,\mathbf{\Theta }}& \leq C\big({\mathbb{P}}_{\mathbf{%
\Theta }}(H>1/8)^{1/p}+\Vert DV\Vert _{q,p,\mathbf{\Theta }}\big) \\
& \leq C\big(\Vert H\Vert _{p,\mathbf{\Theta }}+\Vert H\Vert _{q+1,p(q+1),%
\mathbf{\Theta }}^{q+1}\big)\leq C\Vert H\Vert _{q+1,p(q+1),\mathbf{\Theta }%
}^{q+1}
\end{align*}%
and by using \eqref{est-Hl} one gets
\begin{equation*}
\Vert 1-V\Vert _{q+1,p,\mathbf{\Theta }}\leq C\big({\mathrm{S}}_{G,\mathbf{%
\Theta }}(\bar{p})^{2(q+2)}\mathrm{Q}_{F,G,\mathbf{\Theta }}(q+2,\bar{p}%
)^{8d(q+1)}\Vert F-G\Vert _{q+2,\bar{p},\mathbf{\Theta }}^{q+2}\big)^{q+1}.
\end{equation*}%
But $\Vert F-G\Vert _{q+2,\bar{p},\mathbf{\Theta }}\leq \mathrm{Q}_{F,G,%
\mathbf{\Theta }}(q+2,\bar{p})$, and we get
\begin{align*}
\left\vert \partial _{\alpha }p_{F,\mathbf{\Theta }(1-V)}(y)\right\vert \leq
& C\big({\mathrm{S}}_{F,\mathbf{\Theta }}(p^{\prime })\vee {\mathrm{S}}_{G,%
\mathbf{\Theta }}(p^{\prime })\big)^{a}\mathrm{Q}_{F,G,\mathbf{\Theta }%
}(p^{\prime })^{a}\mathrm{m}_{q+2,p^{\prime }}(\mathbf{\Theta )}^{a}\Vert
F-G\Vert _{q+2,p^{\prime },\mathbf{\Theta }}\times \\
& \times {\mathbb{P}}_{\mathbf{\Theta }}(|F-x|<2)^{b}
\end{align*}%
for $p^{\prime }>d$ and suitable constants $C>0$ and $a>1$ depending on $q,d$. And similarly we get
\begin{align*}
\left\vert \partial _{\alpha }p_{G,\mathbf{\Theta }(1-V)}(y)\right\vert \leq
& C\,{\mathrm{S}}_{G,\mathbf{\Theta }}(p^{\prime })^a\mathrm{Q}_{F,G,\mathbf{%
\Theta }}(q+2,p^{\prime })^a\mathrm{m}_{q+2,p^{\prime }}(\mathbf{\Theta )}%
^{a}\Vert F-G\Vert _{q+2,p^{\prime },\mathbf{\Theta}}\times \\
& \times {\mathbb{P}}_{\mathbf{\Theta }}(|G-x|<2)^{b},
\end{align*}%
with the same constraints for $p^{\prime },C,a$. The statement now follows. $%
\square $

\section{Stochastic equations with jumps}

\label{jumps}

In this section we consider a jump type stochastic differential equation
which has already been considered in \cite{bib:[BCl]}. It is closely related
to piecewise deterministic Markov processes (in fact it is a particular case
of this type of processes). We consider a Poisson point process $p$ with
state space $(E,\mathcal{B}(E)),$ where $E=\mathbb{R}^{d}\times \mathbb{R}%
_{+}.$ We refer to \cite{[IW]} for the notations. We denote by $N$ the
counting measure associated to $p$, we have $N([0,t)\times A)=\#\{0\leq
s<t;p_{s}\in A\}$ for $t\geq 0$ and $A\in \mathcal{B}(E)$. We assume that
the associated intensity measure is given by $\widehat{N}(dt,dz,du)=dt\times
\mu (dz)\times 1_{[0,\infty )}(u)du$ where $(z,u)\in E=\mathbb{R}^{d}\times
\mathbb{R}_{+}$ and $\mu (dz)=h(z)dz.$

We are interested in the solution to the $d$ dimensional stochastic equation
\begin{equation}
X_{t}=x+\int_{0}^{t}\int_{E}c(z,X_{s-})1_{\{u<\gamma
(z,X_{s-})\}}N(ds,dz,du)+\int_{0}^{t}g(X_{s})ds.  \label{eq1}
\end{equation}%
We remark that the infinitesimal generator of the Markov process $X_{t}$ is
given by
\begin{equation*}
L\psi (x)=g(x)\nabla \psi (x)+\int_{\mathbb{R}^{d}}(\psi (x+c(z,x))-\psi
(x))K(x,dz)
\end{equation*}%
where $K(x,dz)=\gamma (z,x)h(z)dz$ depends on the variable $x\in \mathbb{R}%
^{d}.$ See \cite{[F.1]} for the proof of existence and uniqueness of the
solution to (\ref{eq1}).

We describe now our approximation procedure. We consider a non-negative and
smooth function $\varphi :\mathbb{R}^{d}\rightarrow \mathbb{R}_{+}$ such
that $\varphi (z)=0$ for $\left\vert z\right\vert >1$ and $\int_{\mathbb{R}%
^{d}}\varphi (z)dz=1.$ And for $M\in \mathbb{N}$ we denote $\Phi
_{M}(z)=\varphi \ast 1_{B_{M}}$ with $B_{M}=\{z\in \mathbb{R}^{d}:\left\vert
z\right\vert <M\}.$ Then $\Phi _{M}\in C_{b}^{\infty }$ and we have $%
1_{B_{M-1}}\leq \Phi _{M}\leq 1_{B_{M+1}}.$ We denote by $X_{t}^{M}$ the
solution of the equation
\begin{equation}
X_{t}^{M}=x+\int_{0}^{t}\int_{E}c(z,X_{s-}^{M})1_{\{u<\gamma
(z,X_{s-}^{M})\}}\Phi _{M}(z)N(ds,dz,du)+\int_{0}^{t}g(X_{s}^{M})ds.
\label{eq2}
\end{equation}%
In the following we will assume that $\left\vert \gamma (z,x)\right\vert
\leq \overline{C}$ for some constant $\overline{C}.$ Let $%
N_{M}(ds,dz,du):=1_{B_{M+1}}(z)\times 1_{[0,2\overline{C}]}(u)N(ds,dz,du).$
Since $\{u<\gamma (z,X_{s-}^{M})\}\subset \{u<2\overline{C}\}$ and $\Phi
_{M}(z)=0$ for $\left\vert z\right\vert >M+1,$ we may replace $N$ by $N_{M}$
in the above equation and consequently $X_{t}^{M}$ is solution to the
equation
\begin{align*}
X_{t}^{M} &=x+\int_{0}^{t}\int_{E}c_{M}(z,X_{s-}^{M})1_{\{u<\gamma
(z,X_{s-}^{M})\}}N_{M}(ds,dz,du)+\int_{0}^{t}g(X_{s}^{M})ds,\quad \mbox{with}
\\
c_{M}(z,x) &=\Phi _{M}(z)c(z,x).
\end{align*}%
Since the intensity measure $\widehat{N}_{M}$ is finite we may represent the
random measure $N_{M}$ by a compound Poisson process. Let $\lambda _{M}=2%
\overline{C}\times \mu (B_{M+1})=t^{-1}{\mathbb{E}}(N_{M}(t,E))$ and let $%
J_{t}^{M}$ a Poisson process of parameter $\lambda _{M}.$ We denote by $%
T_{k}^{M},k\in \mathbb{N}$ the jump times of $J_{t}^{M}$. We also consider
two sequences of independent random variables $(Z_{k}^{M})_{k\in \mathbb{N}}$
in $\mathbb{R}^{d}$ and $(U_{k})_{k\in \mathbb{N}}$ in $\mathbb{R}_{+}$
which are independent of $J^{M}$ and such that
\begin{equation*}
Z_{k}\sim \frac{1}{\mu (B_{M+1})}1_{B_{M+1}}(z)d\mu (z)\quad \mbox{and}\quad
U_{k}\sim \frac{1}{2\overline{C}}1_{[0,2\overline{C}]}(u)du.
\end{equation*}%
To simplify the notation, we omit the dependence on $M$ for the variables $%
(T_{k}^{M})$ and $(Z_{k}^{M})$. Then equation $(\ref{eq2})$ may be written
as
\begin{equation}
X_{t}^{M}=x+\sum_{k=1}^{J_{t}^{M}}c_{M}(Z_{k},X_{T_{k}-}^{M})1_{(U_{k},%
\infty )}(\gamma (Z_{k},X_{T_{k}-}^{M}))+\int_{0}^{t}g(X_{s}^{M})ds.
\label{eq3}
\end{equation}

In \cite{bib:[BCl]} it is proved that $X_{t}^{M}\rightarrow X_{t}$ in $L^{1}.
$ We study here the convergence in finite variation. Let us give our
hypothesis.

\begin{hypothesis}
\label{3.0} We assume that $\gamma ,g,h$ and $c$ are infinitely
differentiable functions in both variables $z$ and $x$. Moreover we assume
that $g$ and its derivatives are bounded and that $\ln h$ has bounded
derivatives
\end{hypothesis}

\begin{hypothesis}
\label{3.1} We assume that there exist two functions $\overline{\gamma },%
\underline{\gamma }:\mathbb{R}^{d}\rightarrow \mathbb{R}_{+}$ and a constant
$\overline{C}$\ such that
\begin{equation*}
\overline{C}\geq \overline{\gamma }(z)\geq \gamma (z,x)\geq \underline{%
\gamma }(z)\geq 0,\quad \forall x\in \mathbb{R}^{d}
\end{equation*}
\end{hypothesis}

\begin{hypothesis}
\label{3.2}

\begin{itemize}
\item[\textbf{i)}] We assume that there exists a non negative and bounded
function $\overline{c}:\mathbb{R}^{d}\rightarrow \mathbb{R}_{+}$ such that $%
\int_{\mathbb{R}^{d}}\overline{c}(z)d\mu (z)<\infty $ and%
\begin{equation*}
\left\Vert \nabla _{x}c\times (I+\nabla _{x}c)^{-1}(z,x)\right\Vert
+\left\vert c(z,x)\right\vert +\left\vert \partial _{z}^{\beta }\partial
_{x}^{\alpha }c(z,x)\right\vert \leq \overline{c}(z)\quad \forall z,x\in
\mathbb{R}^{d}.
\end{equation*}

\item[\textbf{ii)}] There exists a non negative function $\underline{c}:%
\mathbb{R}^{d}\rightarrow \mathbb{R}_{+}$ such that for every $z\in \mathbb{R%
}^{d}$
\begin{equation*}
\sum_{r=1}^{d}\left\langle \partial _{z_{r}}c(z,x),\xi \right\rangle
^{2}\geq \underline{c}^{2}(z)\left\vert \xi \right\vert ^{2},\quad \forall
\xi \in \mathbb{R}^{d}
\end{equation*}%
and we assume that there exists $\theta >0$ such that
\begin{equation}
\underline{\lim }_{a\rightarrow +\infty }\frac{1}{\ln a}\int_{\{\underline{c}%
^{2}\geq 1/a\}}\underline{\gamma }(z)d\mu (z)=\theta .  \label{eq6}
\end{equation}
\end{itemize}
\end{hypothesis}

\begin{hypothesis}
\label{3.3} We assume that%
\begin{align*}
i) &\qquad \sup_{x,z}\sup_{1\leq |\beta |\leq l}|\partial _{\beta ,z}\ln
\gamma (z,x)| <\infty , \\
ii)&\qquad \sup_{z^{\ast }\in \mathbb{R}^{d}}\int_{B(z^{\ast },1)}\overline{%
\gamma }(z)d\mu (z) <+\infty , \\
iii)&\qquad \int_{\mathbb{R}^{d}}\overline{\gamma }_{\ln }^{x,l}(z)\overline{%
\gamma }(z)d\mu (z) <\infty
\end{align*}%
with $\overline{\gamma }_{\ln }^{x,l}(z)=\sup_{x}\sup_{1\leq |\beta |\leq
l}|\partial _{\beta ,x}\ln \gamma (z,x)|.$
\end{hypothesis}

We are now able to give our convergence result.

\begin{theorem}
\label{Convergence} Suppose that Hypothesis \ref{3.0}-\ref{3.3} hold. Then
for every $t>0$ 
one has
\begin{equation*}
\lim_{M\to \infty}d_{TV}(X_{t},X_{t}^{M})=0.
\end{equation*}
\end{theorem}

\textbf{Proof}. The proof is an easy consequence of the results from \cite%
{bib:[BCl]}, we use the estimates obtained there.

\smallskip

\textbf{Step 1.} In \cite{bib:[BCl]} Lemma 4 one proves that $%
X_{t}^{M}\rightarrow X_{t} $ in $L^{1}$ and then $\lim_{M\rightarrow \infty
}d_{1}(X_{t},X_{t}^{M})=0.$

\smallskip

\textbf{Step 2.} Following \cite{bib:[BCl]}, we consider an alternative
representation of the law of $X_{t}^{M}$. The random variable $X_{t}^{M}$
solution to (\ref{eq3}) is a function of $(Z_{1}\ldots ,Z_{J_{t}^{M}})$ but
it is not a simple functional, as defined in Section \ref{th-abstract},
because the coefficient $c_{M}(z,x)1_{(u,\infty )}(\gamma (z,x))$ is not
differentiable with respect to $z$. In order to avoid this difficulty we use
the following alternative representation. Let $z_{M}^{\ast }\in \mathbb{R}%
^{d}$ such that $\left\vert z_{M}^{\ast }\right\vert =M+3$. We define
\begin{equation*}
\begin{array}{rcl}
q_{M}(z,x) & := & \displaystyle\varphi (z-z_{M}^{\ast })\theta _{M,\gamma
}(x)+\frac{1}{2\overline{C}\mu (B_{M+1})}1_{B_{M+1}}(z)\gamma (z,x)h(z)
\smallskip \\
\theta _{M,\gamma }(x) & := & \displaystyle\frac{1}{\mu (B_{M+1})}%
\int_{\{\left\vert z\right\vert \leq M+1\}}\Big(1-\frac{1}{2\overline{C}}%
\gamma (z,x)\Big)\mu (dz).%
\end{array}
\end{equation*}
We recall that $\varphi $ is a non-negative and smooth function with $\int
\varphi =1$ and which is null outside the unit ball. Moreover since, $0\leq
\gamma (z,x)\leq \overline{C}$ and then $1\geq \theta _{M,\gamma }(x)\geq 1/2
$. By construction the function $q_{M}$ satisfies $\int q_{M}(x,z)dz=1.$
Hence we can easily check (see \cite{bib:[BCl]} for a complete proof) that
\begin{equation}
{\mathbb{E}}(f(X_{T_{k}}^{M})\mid X_{T_{k}-}^{M}=x)=\int_{{\mathbb{R}}%
^{d}}f(x+c_{M}(z,x))q_{M}(z,x)dz.  \label{condXM}
\end{equation}%
From the relation (\ref{condXM}) we construct a process $(\overline{X}%
_{t}^{M})$ equal in law to $(X_{t}^{M})$ in the following way. We denote by $%
\Psi _{t}(x)$ the solution of $\Psi _{t}(x)=x+\int_{0}^{t}g(\Psi _{s}(x))ds.$
We assume that the times $T_{k},k\in \mathbb{N}$ are fixed and we consider a
sequence $(z_{k})_{k\in \mathbb{N}}$ with $z_{k}\in \mathbb{R}^{d}.$ Then we
define $x_{t},t\geq 0$ by $x_{0}=x$ and, if $x_{T_{k}}$ is given, then
\begin{eqnarray*}
x_{t} &=&\Psi _{t-T_{k}}(x_{T_{k}})\quad T_{k}\leq t<T_{k+1}, \\
x_{T_{k+1}} &=&x_{T_{k+1}^{-}}+c_{M}(z_{k+1},x_{T_{k+1}^{-}}).
\end{eqnarray*}%
We remark that for $T_{k}\leq t<T_{k+1},x_{t}$ is a function of $%
z_{1},...,z_{k}.$ Notice also that $x_{t}$ solves the equation
\begin{equation*}
x_{t}=x+\sum_{k=1}^{J_{t}^{M}}c_{M}(z_{k},x_{T_{k}^{-}})+%
\int_{0}^{t}g(x_{s})ds.
\end{equation*}%
We consider now a sequence of random variables $(\overline{Z}_{k}),k\in
\mathbb{N}^{\ast }$ and we denote $\mathcal{G}_{k}=\sigma (T_{p},p\in
\mathbb{N})\vee \sigma (\overline{Z}_{p},p\leq k)$ and $\overline{X}%
_{t}^{M}=x_{t}(\overline{Z}_{1},...,\overline{Z}_{J_{t}^{M}}).$ We assume
that the law of $\overline{Z}_{k+1}$ conditionally on $\mathcal{G}_{k}$ is
given by
\begin{equation*}
P(\overline{Z}_{k+1}\in dz\mid \mathcal{G}_{k})=q_{M}(x_{T_{k+1}^{-}}(%
\overline{Z}_{1},...,\overline{Z}_{k}),z)dz=q_{M}(\overline{X}%
_{T_{k+1}^{-}}^{M},z)dz.
\end{equation*}%
Clearly $\overline{X}_{t}^{M}$ satisfies the equation
\begin{equation}
\overline{X}_{t}^{M}=x+\sum_{k=1}^{J_{t}^{M}}c_{M}(\overline{Z}_{k},%
\overline{X}_{T_{k}-}^{M})+\int_{0}^{t}g(\overline{X}_{s}^{M})ds.
\label{eq4}
\end{equation}%
Notice that $\overline{X}_{t}^{M}$\ is a piecewise deterministic Markov
process, but not a completely general one because the intensity of the law
of the jump times $T_{k},k\in {\mathbb{N}}$ does not depend on the position
of the particle $\overline{X}_{t}^{M}.$ We think that the more general case
may also be considered using similar arguments but we leave this out here.

\smallskip

\textbf{Step 3.} We will use the integration by parts formulae from Section %
\ref{th-abstract} with the random variable $V=(V_{1},...,V_{J})$ replaced by
$(\overline{Z}_{1},...,\overline{Z}_{J_{t}^{M}})$ with fixed $M$ and $t>0.$
We use the weight $\pi _{k,r}=\Phi _{M}(\overline{Z}_{k}),k\in {\mathbb{N}}%
,r=1,...,d$ and the Malliavin derivative is
\begin{equation*}
D_{k,r}=\pi _{k,r}\partial _{\overline{Z}_{k}^{r}}.
\end{equation*}%
In fact we will work conditionally to the time grid $T_{k},k\in {\mathbb{N}}$
but all the constants coming on are independent of the time grid (as well as
on $M$ and $t>0$) so we do not mention this in the notation.

We will use several estimates obtained in \cite{bib:[BCl]}. First, by Lemma
7 and Lemma 13 in \cite{bib:[BCl]}, for every $p\geq 1,l\in {\mathbb{N}}$ we
have
\begin{equation}
{\mathbb{E}}(\vert \overline{X}_{t}^{M}\vert _{l}^{p})+{\mathbb{E}}(\vert L%
\overline{X}_{t}^{M}\vert _{l}^{p})\leq C_{l,p}  \label{eq5}
\end{equation}
with $C_{l,p}\in \mathcal{C}(d,l,p).$ So hypothesis (\ref{New3}) holds for
every $r.$

We discuss now the non degeneracy property. We consider the tangent flow $%
Y_{t}^{M}$ solution to
\begin{equation}
Y_{t}^{M}=I+\sum_{k=1}^{J_{t}^{M}}\nabla _{x}c_{M}(\overline{Z}_{k},%
\overline{X}_{T_{k}-}^{M})Y_{T_{k}-}^{M}+\int_{0}^{t}\nabla _{x}g(\overline{X%
}_{s}^{M})Y_{s}^{M}ds.  \label{flot}
\end{equation}%
Since $\left\Vert \nabla _{x}c\times (I+\nabla _{x}c)^{-1}(z,x)\right\Vert
\leq \overline{c}(z)$ it follows that $Y_{t}^{M}$ is invertible; we denote
by $\widehat{Y}_{t}^{M}$ its inverse. Then it is proved in \cite{bib:[BCl]}
that $D_{k,r}\overline{X}_{t}^{M,r^{\prime }}=\pi _{k}(Y_{t}^{M}\nabla
_{z}c_{M}(\overline{Z}_{k},\overline{X}_{T_{k}^{-}}^{M}))_{r^{\prime },r}$
and moreover, if $\lambda _{t}$ denotes the lower eigenvalue of the
Malliavin covariance matrix of $\overline{X}_{t}^{M}$ we have
\begin{equation*}
\rho _{t}^{M}\geq \Vert \widehat{Y}_{t}^{M}\Vert
^{-2}\sum_{k=1}^{J_{t}^{M}}1_{B_{M-1}}(\overline{Z}_{k})\underline{c}^{2}(%
\overline{Z}_{k}).
\end{equation*}%
Then
\begin{align*}
{\mathbb{P}}(\sigma _{\overline{X}_{t}^{M}} \leq \varepsilon ) &\leq {%
\mathbb{P}}\big(\Vert \widehat{Y}_{t}^{M}\Vert
^{-2}\sum_{k=1}^{J_{t}^{M}}1_{B_{M-1}}(\overline{Z}_{k})\underline{c}^{2}(%
\overline{Z}_{k})\leq \varepsilon ^{1/d}\big) \\
&\leq {\mathbb{P}}\Big(\sum_{k=1}^{J_{t}^{M}}1_{B_{M-1}}(\overline{Z}_{k})%
\underline{c}^{2}(\overline{Z}_{k})\leq \varepsilon ^{1/2d}\Big)+ {\mathbb{P}%
}\big(\Vert \widehat{Y}_{t}^{M}\Vert ^{-2}\leq \varepsilon ^{1/2d}\big) \\
&\leq {\mathbb{P}}\Big(\sum_{k=1}^{J_{t}^{M}}\Phi _{M}(\overline{Z}_{k})%
\underline{c}^{2}(\overline{Z}_{k})\leq \varepsilon ^{1/2d}\Big)+{\mathbb{P}}%
\big(\Vert Y_{t}^{M}\Vert ^{2}\geq \varepsilon ^{1/2d}\big).
\end{align*}%
In \cite{bib:[BCl]} one proves that $\sup_{M}{\mathbb{E}}\Vert
Y_{t}^{M}\Vert ^{2p}<\infty $ for every $p\geq 1$ so that
\begin{equation*}
\limsup _{\varepsilon \rightarrow 0}\limsup _{M\rightarrow 0}P(\left\Vert
Y_{t}^{M}\right\Vert ^{2}\geq \varepsilon ^{1/2d})=0.
\end{equation*}
One also proves in Lemma 5 from \cite{bib:[BCl]} that $%
\sum_{k=1}^{J_{t}^{M}}\Phi _{M}(\overline{Z}_{k})\underline{c}^{2}(\overline{%
Z}_{k})$ has the same law as
\begin{equation*}
\sum_{k=1}^{J_{t}^{M}}\Phi _{M}(Z_{k})\underline{c}^{2}(Z_{k})1_{[0,U_{k}]} %
\big(\vert \gamma (Z_{k},X_{T_{k}-}^{M})\vert \big),
\end{equation*}
so
\begin{equation*}
{\mathbb{P}}\Big(\sum_{k=1}^{J_{t}^{M}}\Phi _{M}(\overline{Z}_{k})\underline{%
c}^{2}(\overline{Z}_{k})\leq \varepsilon ^{1/2d}\Big) ={\mathbb{P}}\Big(%
\sum_{k=1}^{J_{t}^{M}}\Phi _{M}(Z_{k})\underline{c}^{2}(Z_{k})1_{[0,U_{k}]}(%
\left\vert \gamma (Z_{k},X_{T_{k}-}^{M})\right\vert )\leq \varepsilon ^{1/2d}%
\Big).
\end{equation*}%
Let us denote
\begin{equation*}
N_{M}(t)=\sum_{k=1}^{J_{t}^{M}}\Phi _{M}(Z_{k})\underline{c}%
^{2}(Z_{k})1_{[0,U_{k}]}(\left\vert \gamma (Z_{k},X_{T_{k}-}^{M})\right\vert
)\quad\mbox{and}\quad U_{M}(t)=t\int_{B_{M-1}^{c}}\underline{c}^{2}(z)%
\underline{\gamma }(z)d\mu (z).
\end{equation*}
In the final part of the proof of Lemma 16 in \cite{bib:[BCl]} one shows
that if $p/t<\theta $ (with $\theta $ from (\ref{eq6})) then ${\mathbb{E}}%
((N_{M}(t)+U_{M}(t))^{-p})\leq C_{p}.$ Since $\lim_{M\rightarrow \infty
}U_{M}(t)=0$ one has for every fixed $\varepsilon >0$
\begin{equation*}
\limsup_{M\rightarrow \infty }{\mathbb{P}}(N_{M}(t)<\varepsilon )
=\limsup_{M\rightarrow \infty }{\mathbb{P}}(N_{M}(t)+U_{M}(t)<\varepsilon )
\leq \varepsilon ^{p}{\mathbb{E}}((N_{M}(t)+U_{M}(t))^{-p})\leq
C_{p}\varepsilon ^{p}.
\end{equation*}%
So if $\theta >0$ we take $p<\theta t$ and we obtain
\begin{equation*}
\limsup_{\varepsilon \rightarrow 0} \limsup_{M\rightarrow \infty }{\mathbb{P}%
}(N_{M}(t)<\varepsilon )=0
\end{equation*}%
so that hypothesis (\ref{New4}) is also verified. Now the conclusion follows
from Theorem \ref{General}. $\square $

\end{document}